\documentclass[aap,MSNbibl]{arximspdf}
\usepackage{graphicx}
\usepackage{upgreek}

\doi{10.1214/13-AAP950} 
\volume{24}
\issue{3}
\pubyear{2014}
\firstpage{1298}
\lastpage{1328}

\makeatletter
\DeclareMathAlphabet{\mathbbl}{U}{mt2hrb}{m}{n}
\SetMathAlphabet{\mathbbl}{bold}{U}{mt2hrb}{b}{n}
\newcommand{\mathbbe}{\mathbbl{e}}
\newcommand{\rrVert}{\Vert}
\newcommand{\rrvert}{\vert}
\newcommand{\llVert}{\Vert}
\newcommand{\llvert}{\vert}
\newproclaim{remark}{Remark}
\newproclaim{Imremark}{Important Remark}
\newproclaim{example}{Example}
\newtheorem{theorem}{Theorem}
\newtheorem{corollary}[theorem]{Corollary}
\newtheorem{proposition}[theorem]{Proposition}

\newproclaim{definition}{Definition}
\makeatother

\begin{document}
\begin{frontmatter}

\title{Propagation of chaos in neural fields}
\runtitle{Stochastic neural fields theory}
\pdftitle{Propagation of chaos in neural fields}

\begin{aug}
\author{\fnms{Jonathan} \snm{Touboul}\corref{}\ead[label=e1]{jonathan.touboul@college-de-france.fr}\ead[label=u1,url]{http://mathematical-neuroscience.net/team/jonathan}}
\runauthor{J. Touboul}
\affiliation{CIRB - Coll\`ege de France and INRIA Paris-Rocquencourt}
\address{The  Mathematical Neuroscience Team\\
CIRB and INRIA\\
Coll\`ege de France\\
11 Place Marcelin Berthelot\\
75005 Paris\\
France\\
\printead{e1}\\
\printead{u1}} 
\end{aug}

\received{\smonth{7} \syear{2012}}
\revised{\smonth{3} \syear{2013}}

%
\begin{abstract}
We consider the problem of the limit of bio-inspired spatially extended
neuronal networks including an infinite number of neuronal types (space
locations), with space-dependent propagation delays modeling neural
fields. The propagation of chaos property is proved in this setting
under mild assumptions on the neuronal dynamics, valid for most models
used in neuroscience, in a mesoscopic limit, the \textit{neural-field
limit}, in which we can resolve the quite fine structure of the
neuron's activity in space and where averaging effects occur. The
mean-field equations obtained are of a new type: they take the form of
well-posed infinite-dimensional delayed integro-differential equations
with a nonlocal mean-field term and a singular spatio-temporal Brownian
motion. We also show how these intricate equations can be used in
practice to uncover mathematically the precise mesoscopic dynamics of
the neural field in a particular model where the mean-field equations
exactly reduce to deterministic nonlinear delayed integro-differential
equations. These results have several theoretical implications in
neuroscience we review in the discussion.
\end{abstract}

%
\pdfkeywords{60F99, 60B10, 34C15, Mean-field limits, propagation of chaos, delayed stochastic differential equations,
infinite-dimensional stochastic processes, neural fields}

\end{frontmatter}

{\bf Update:} This document was updated in May 2016, to clarify the case of an extension suggestion in Appendix B that was erroneous in its original form, and to provide details on the estimates used in the proof of Theorem~\ref{thmExistenceUniquenessSpace}. We thank Wilhelm Stannat and Fran\c cois Delarue for their remarks.

\section*{Introduction}\label{sec1}

The brain's activity is the result of the complex interplay of
different cells, in particular neurons, electrical cells that manifest
highly complex nonlinear behaviors characterized by the intense
presence of noise. Neurons form large population assemblies at the
scale of which emerge reliable and adapted responses to stimuli. Such
local neural populations, often termed cortical columns, have a
diameter of about 50~$\upmu$m to 1~mm, contain a few thousand to one
hundred thousand neurons and are in charge of specific
functions~\cite{mountcastle97}. The interaction of several columns at
different spatial locations allows processing of the complex sensory or
cortical information and supports brain function. Such groups of
cortical columns organize on the surface of the cortex and form
spatially extended structures called \emph{neural fields}, the activity
of which is precisely at the scale most usual imaging techniques (e.g.,
EEG/MEG, optical imaging) record relevant phenomena, and also
correspond to anatomical information revealed experimentally. A
paradigmatic example is given by the primary visual cortex of certain
mammals. In such cortical areas, neurons organize into columns
responding preferentially to specific orientations in visual stimuli
and display specific connection patterns
\cite{hubel-wiesel-etal78,bosking-zhang-etal97}. The communication
between neurons is characterized by a delay due to the transport of
information through axons and to the typical time the synaptic
machinery needs to transmit it. These delays have a clear role in
shaping the neuronal activity, as established by different authors;
see, for example,~\cite{coombes-laing11,series-fregnac02}. In such
structures, several highly populated columns interact, and the number
of neurons in each column is orders of magnitude higher than the number
of columns (e.g., orientations) involved. A variety of important brain
states rely on the coordinated behaviors of large neural assemblies and
recently raised the interest of physiologists and computational
neuroscientists. Among these, we shall cite the rapid complex answers
to specific stimuli~\cite{thorpe-delorme-etal01}, decorrelated
activity~\cite{ecker-berens-etal10,renart-de-la-rocha-etal10}, large
scale oscillations~\cite{buszaki06},
synchronization~\cite{izhikevich-polychronization06} and
spatio-temporal pattern
formation~\cite{ermentrout-cowan80,coombes-owen05}.

The mathematical and computational analysis of the dynamics of neural
fields relies almost exclusively on the use of heuristic models since
the seminal work of Wilson and Cowan \cite{wilson-cowan73} and Amari
\cite{amari72}. This approach implicitly considers that
averaging effects counterbalance the prominent noisy aspect of in vivo
firing observed experimentally, and describes the mesoscopic cortical
activity through a deterministic, scalar variable whose dynamics is
given by integro-differential equations. This model was widely studied
analytically and numerically, and successfully accounted for
hallucination patterns, binocular rivalry and
synchronization~\cite{laing-troy-etal02,ermentrout-cowan79}. Justifying
these models starting from biologically realistic settings has since
then been a great endeavor~\cite{bressloff12}.

In this manuscript we undertake a rigorous analysis of neural fields.
From the biological viewpoint, these are spatially extended cortical
structures made of several highly populated neuronal ensembles (the
neural \emph{populations}) in charge of specific functions. From the
mathematical viewpoint, neural fields are adequately described as the
limit of a set of nonlinear interacting stochastic processes (generally
governing the neuron's electrical potential and related variables)
gathering into different homogeneous populations at specific locations
on the cortex. Neurons in each population have similar dynamics and
communicate with neurons of different populations depending on the
respective positions of the populations on the cortex and after a
specific time delay. In what we will call the \emph{neural field
limit}, both the number of neurons and the number of populations tend
to infinity so that the populations completely cover a continuous space
(a piece of cortex or a functional space).

This problem is evocative of statistical fluid mechanics and more
generally interacting particle systems, and as such has been widely
studied in mathematics and physics, chiefly motivated by thermodynamics
or fluid dynamics questions. In particular, the probability
distribution of a typical set of particles in the limit where the total
number of particles goes to infinity, and fluctuations around this
limit where characterized for a number of
models~\cite{mckean66,dobrushin70,tanaka83,sznitman89,sznitman84a}. It
was shown in several contexts that when considering that all particles
have independent identically distributed initial conditions
(\emph{chaotic} initial conditions), then in the limit where the number
of particles tends to infinity, the behavior of a few particles remains
independent as time goes by, and all particles have the same
probability distribution, which is the solution of a nonlinear Markov
equation, often referred to as the \emph{McKean--Vlasov} equation. The
underlying biological problem motivates the \hyperref[sec1]{Introduction} of a notion of
spatial labeling of the (fixed) neurons, involving two mathematical
aspects that were not covered in the literature. First is the fact that
this induces the presence of infinitely many types of neurons
(corresponding to the column neurons belong to), and second is the fact
that since neurons communicate through the emission of electrical
impulses transported at finite speed through the axons, space-dependent
delays occur in the communication between two cells. These two aspects
necessitate the development of the propagation of chaos theory toward
infinite-dimensional functional settings that we aim at achieving in
the present manuscript. We will show that in the neural field limit,
the propagation of chaos property holds. Moreover, the activity is
shown to converge in a certain sense toward the solution of a new
object, a delayed integro-differential mean-field equation with
space-dependent delays. This object is substantially different from the
usual McKean--Vlasov limits: beyond the presence of delays, the neural
field limit regime is at a mesoscopic scale where averaging effects
locally occur, but is fine enough to resolve brain's structure and its
activity, resulting in the presence of an integral term over space. The
speed of convergence toward the mean-field equations is quantified and
involves two terms, one governing the averaging effect in each
population and the second corresponding to the continuum limit. In the
neural field regime, the limit equations are very singular; in
particular, trajectories are not measurable with respect to the space.
These limits are very hard to analyze at this level of generality.
However, in the type of models usually considered in the study of
neural fields, namely the firing-rate model, we show in a companion
article~\cite{touboulNeuralFieldsDynamics11} that the behavior can be
rigorously and exactly reduced to a system of deterministic
integro-differential equations that are compatible with the usual
Wilson and Cowan system in the zero noise limit. Noise intervenes in
these equations a nonlinear fashion, fundamentally shaping in the
macroscopic dynamics.

The paper is organized as follows. We start in
Section~\ref{secMathSetting} by describing the mathematical setting of
the study, abstracting classical relevant neuronal models that are
specified and reviewed in Appendix~\ref{secMotiv}, and more general
models are considered in Appendix~\ref{appendStoSynapses}. We then
analyze the integro-differential delayed McKean--Vlasov equations that
will constitute our limit neural field equation in
Section~\ref{secExistenceUniquenessSpace} and demonstrate in particular
their well-posedness, before addressing in
Section~\ref{secPropaChaSpace} the propagation of chaos property and
convergence of the network equations toward the solutions of the
mean-field equation. In Section~\ref{secFiringRates} we illustrate how
this approach can be used in practice to analyze the effect of the
parameters on the dynamics of the system in a particular example,
reviewing some results of~\cite{touboulNeuralFieldsDynamics11} afresh
on a new example where noise, delays and spatial structure interact to
shape the mesoscopic response of the neural field. The results of the
mathematical analysis are then confronted to different recent
experimental observations on collective dynamics of neural fields in
the brain, and a few open problems of interest are discussed in the
conclusion Section~\ref{secdiscussion}.

%
\begin{figure}

\includegraphics{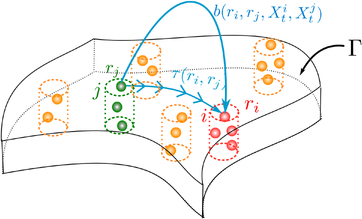}

\caption{A typical architecture of neural fields: cylinders represent
neural populations as cortical columns spanning across the cortex.
Neuron~$i$ (red population at $r_i\in\Gamma$) receives a spike from
neuron~$j$ (green population at $r_{j}\in\Gamma$) after a delay $\tau
(r_i,r_j)$ creating a current $b(r_i,r_j,X^i_t,X^j_t)$.} \label{figNeurons}
\end{figure}

\section{Mathematical setting}\label{secMathSetting}

Throughout the manuscript, we work in a complete probability space
$(\Omega, \mathcal{F},\mathbb{P})$ endowed with a filtration $
(\mathcal{F}_t )_t$ satisfying the usual conditions. We consider a
spatially extended network (see Figure~\ref{figNeurons}) composed of $N$ neurons, each neuron
belonging to one of $P(N)$ populations characterized by their locations
$(r_1,\ldots, r_{P(N)}) \in\Gamma^{P(N)}$ on the cortex (or the feature
space) $\Gamma$, a finite-dimensional compact set.\footnote{When
considering $\Gamma$ as the cortex, it will be a compact subset of
$\mathbb{R}^{q}$, $q=2$ or $3$, and when considering that populations
are defined by the neuron's function, the shape of $\Gamma$ can take
different forms depending on the geometry of the feature space. For
instance, in the case of the primary visual area, neurons code for a
preferred orientation of a visual stimuli that can be represented in
the torus $\Gamma=\mathbb{S}^1$.} The state of each neuron $i$ in the
network is described by a $d$-dimensional variable $X^i\in
E:=\mathbb{R}^d$, typically corresponding to the membrane potential of
the neuron and possibly additional variables such as those related to
ionic concentrations and gated channels described in
Appendix~\ref{secMotiv}, and satisfy the network equations
%
\begin{eqnarray}
\label{eqNetworkSpace} d X^{i,N}_t &=& \Biggl(f\bigl(r_{\alpha},t,X^{i,N}_t
\bigr)\nonumber
\\
&&\hspace*{6pt}
{} + \frac{{1}}{P(N)}\sum_{\gamma=1}^{P(N)}\sum_{p(j)=\gamma} \frac
{1}{N_{\gamma}} b\bigl(r_{\alpha},r_{\gamma},X^{i,N}_t,X^{j,N}_{t-\tau
(r_{\alpha},r_{\gamma})}
\bigr) \Biggr) \,dt
\\
&&{}+ \sigma(r) \,dW^{i}_t,\nonumber
\end{eqnarray}
where $f(r,t,x)\dvtx \Gamma\times\mathbb{R}\times E\mapsto E$ governs
the intrinsic dynamics of each cell, $(W_t^i)$~is a sequence of
$m$-dimensional Brownian motions modeling the external noise and
$\sigma(r)\dvtx \Gamma\mapsto\mathbb{R}^{d\times m}$ a bounded and
measurable function of $r\in\Gamma$ modeling the level of noise at each
space location and $b(r,r',x,y)\dvtx \Gamma^2\times E^2\mapsto E$ the
interaction function of a neuron located at $r'$ with voltage $y$ on a
neuron at location $r$ with voltage $x$. The function $\tau(r,r')\dvtx
\Gamma ^2\mapsto \mathbb{R}^+$ is the interaction delay between neurons
located at $r$ and those at $r'$ which is assumed to be a regular
function of its two variables. We assume that all delays are bounded by
a finite quantity $\tau$. The quantity $r_{\alpha}$ is called the
\emph{location} of the population, and $\alpha$ is the population
label. For a neuron $i$ in the network, the population function $p\dvtx
\mathbb{N}\mapsto\mathbb{N}$ associates to a neuron $i$ the population
$\alpha$ it belongs to. The number of neurons in each population in a
network of size $N$ defines a sequence of population size $(N_1(N),
\ldots, N_{P(N)}(N))$ [we hence have $\sum_{\gamma=1}^{P(N)}
N_{\gamma}(N)=N$] corresponding to the number of neurons in population
$\gamma$ when the network size is equal to $N$. The number of
populations $P(N)$ and the number of neurons in each of these
populations is assumed to be deterministic.\footnote{It is easy to
generalize to random population number and population size.} The
interaction term presents a scaling factor $\frac{1} {P(N)N_{\gamma}}$
ensuring the boundedness of the input received by neurons from
population $\gamma$ to the other populations, a biological fact related
to the brain function and to the finiteness of the resources available
for the synaptic transmission.

The different locations $r_{\gamma}$ of the populations are related to
the organization of the neurons on the space $\Gamma$. These locations
are distributed according to a~specific probability measure\vspace*{1pt} $\lambda$
on $\Gamma$.\footnote{In the example of the visual area V1, $\lambda$ is
the uniform measure on $\mathbb{S}^1$.} The locations of the $P(N)$
populations, $(r_1,\ldots, r_{P(N)}) \in\Gamma^{P(N)}$, are assumed
to be randomly and independently drawn in $\Gamma$ according to the
probability $\lambda(dr)$ in a different probability space
$(\Omega',\mathcal{F}',\mathbb{P}')$. We will denote by $\mathcal
{E}$ the expectation
over the realizations of the space locations~$(r_{\alpha})$.

It is clear that the larger the number of populations, the smaller the
mean number of neurons per populations. The number of populations will
hence compete with the typical number of neurons per population and
hence with averaging effects. In the present article, motivated by the
fact that the number of neurons in each population is orders of
magnitude larger than the number of populations (see,
e.g.,~\cite{changeux06}), we will make the following assumption,
referred to as the \emph{neural field limit}:
\begin{equation}
\label{eqPopulationEstimate}
\mathbbe (N):=\frac{1} {P(N)} \sum_{\gamma=1}^{P(N)}
\frac{1} {
N_{\gamma}(N)} \mathop{\longrightarrow}\limits
_{N\to\infty} 0.
\end{equation}

In the case of an infinite number of populations, this assumption
ensures heuristically most populations are made of a diverging number
of neurons.\footnote{If all populations have approximately the same
number of neurons, each $N_k(N)$ will be of the order $N/P(N)$, and the
condition~(\ref{eqPopulationEstimate}) is satisfied when $P(N)=o(N)$.
The condition also ensures the size of most populations tend to
infinity. Indeed, for instance, if all but one population contains just
$1$ neuron, the last population contains $N-P(N)$ neurons, and the sum
is equal to $1-1/P(N) + 1/(P(N)(N-P(N))) \geq1-1/P(N)$ which will not
tend to zero.}

The parameters of the system are assumed to satisfy the following
assumptions:\vadjust{\goodbreak}
\begin{longlist}[(H5)]
\item[(H1)]\label{AssumpLocLipschSpace} $f(r,t,\cdot)$ is uniformly $K_f$
    Lipschitz-continuous.

\item[(H2)]\label{AssumpLocLipschbSpace} $b(r,r',\cdot,\cdot)$ is
    uniformly $L$-Lipschitz-continuous.

\item[(H3)]\label{AssumpbBoundSpace} There exists a $\tilde{K}>0$ such
    that
\[
\bigl| b\bigl(r,r',x,z\bigr)\bigr|^2 \leq\tilde{K}
\bigl(1+| x |^2\bigr).
\]
\item[(H4)]\label{AssumpLinearGrowth} The drift satisfies uniformly in
    space ($r$) and time ($t$), the inequality
\[
\bigl| f(r,t,x)\bigr|^2 \leq C \bigl(1+| x |^2
\bigr).
\]

\item[(H5)]\label{AssumpSpaceContinuity} The drift, delay, diffusion and
    coupling functions are regular with respect to space variables
    $(r,r')\in\Gamma^2$ (at least measurable, in practice generally
    assumed continuous).
\end{longlist}

Let us first state the following proposition ensuring well-posedness of
the network system under the assumptions of the section:

\begin{proposition}\label{proExistenceUniquenessNetwork}
Let $(X^0_t)_{t\in[-\tau,0]}$ a square integrable process with values
in~$E^N$. Under the assumptions of the section, there exists a unique
strong solution to the network equations (\ref{eqNetworkSpace}) with
initial condition $X^0$, which is square integrable and defined for all
times.
\end{proposition}

The proof of this proposition is a direct application of Da
Prato~\cite{da-prato92} as used by Mao~\cite{mao02}, and essentially
uses the same arguments as those\vadjust{\goodbreak} of the proof
Theorem~\ref{thmExistenceUniquenessSpace}. The interested reader is
invited to follow the steps of the demonstration of that theorem to
prove Proposition~\ref{proExistenceUniquenessNetwork}.

We are interested in the limit of such systems as the number of neurons
$N$ goes to infinity, under the neural field limit condition.

Let us start by briefly bring some results from the analysis of finite
populations networks [i.e., the case where $P(N)$ remains finite as
$N\to\infty$], which can be seen as a particular case of the current
setting under the assumption that $\lambda$ is a sum of Dirac masses.
In that case, the neural field regime~(\ref{eqPopulationEstimate})
amounts assuming that the number of neurons in each population tends to
infinity. Standard theory proves that the network converges toward $P$
coupled McKean--Vlasov equations
\begin{eqnarray*}
d \bar{X}_t(r_{\alpha})&=&f\bigl(r,t,\bar{X}_t(r_{\alpha})
\bigr) \,dt + \sigma(r_{\alpha}) \,dW_t^{\alpha}
\\
&&{}+ \frac{1} P \sum_{\gamma=1}^P
\mathbb{E}_{\bar{Z}}\bigl[b\bigl(r_{\alpha
},r_{\gamma},
\bar{X}_t(r_{\alpha}),\bar{Z}_{t-\tau(r_{\alpha
},r_{\gamma})}(r_{\gamma})
\bigr)\bigr] \,dt,
\end{eqnarray*}
where $(W^{\alpha}_t)$ are $P$ independent Brownian motions. This model
can be seen as a~discrete approximation of the continuous neural field.
When the asymptotic number of populations is infinite, corresponding
heuristically to refining the spatial discretization (or increasing the
number of populations), one is likely to face two main difficulties:
(i)~the network equations will involve an infinite number of
independent Brownian motions, one for each space location, and (ii)~it
will involve a~limit, as $P$ goes to infinity, of a~sum of the
mean-field interaction terms [it is, rather, a~simultaneous limit under
the scaling property~(\ref{eqPopulationEstimate})].

\begin{remark*}
Note that the infinite number of independent Brownian motions is not a~technical artifact, but a~fact related to the very nature of the
problem: distinct neurons are driven by independent Brownian motions
whatever their respective locations on the neural field $\Gamma$, and
no spatial continuity or measurability is to be expected in the
solution of the limit equations.
\end{remark*}

In order to handle the first point, we introduce a~particular object,
the \emph{spatially chaotic}\footnote{We use the term \emph{chaotic} in
the statistical physics sense as understood by Boltzmann's in his
notion of molecular chaos ``Sto\ss zahlansatz.''} Brownian
motion on $\Gamma$, a~two-parameter process
$(t,r)\in\mathbb{R}^+\times\Gamma \mapsto W_t(r)$ such that for any
fixed $r\in\Gamma$, the process $t\mapsto W_t(r)$ is a~$d$-dimensional
standard Brownian motion, and for $r\neq r'$ in $\Gamma$, the processes
$W_t(r)$ and $W_t(r')$ are independent. This process is relatively
singular seen as a~spatio-temporal process: in particular, it is not
measurable with respect to the Borel algebra $\mathcal{B}(\Gamma)$
of~$\Gamma$. This object, defined as a~collection of independent Brownian
motions, clearly exists. More generally, in what follows,\vadjust{\goodbreak} a~process
$\zeta_t(r)$ will be termed \emph{spatially chaotic} if the processes
$\zeta_t(r)$ and $\zeta_t(r')$ are independent for any $r\neq r'$.

We will show that the network equations~(\ref{eqNetworkSpace}) satisfy
the propagation of chaos property in the limit where $N$ goes to
infinity under the neural field assumption, and that the state of the
network converges toward a~very particular McKean--Vlasov equation
involving a~spatially chaotic Brownian motion. In detail, for almost
all realizations of the spatial locations $(r_{\gamma},\gamma\in
\mathbb{N})$ i.i.d. with law $\lambda$, the asymptotic law of neurons
located at $r$ in the support of $\lambda$ will be measurable with
respect to $(\Gamma,\mathcal{B}(\Gamma))$ and converge toward the
stochastic neural field mean-field equation with delays
\begin{eqnarray}
\label{eqMFESpace} d \bar{X}_t(r)&=&f\bigl(r,t,\bar{X}_t(r)
\bigr) \,dt + \sigma(r) \,dW_t(r)
\nonumber\\[-8pt]\\[-8pt]
&&{} + \int_{\Gamma}
\mathbb{E}_{\bar{Z}}\bigl[b\bigl(r,r',\bar{X}_t(r),
\bar {Z}_{t-\tau(r,r')}\bigl(r'\bigr)\bigr)\bigr] \,d\lambda
\bigl(r'\bigr) \,dt,\nonumber
\end{eqnarray}
where $(W_t(r))_{t\geq0, r\in\Gamma}$ is a~spatially chaotic Brownian,
and the process $(\bar{Z})$ is independent and has the same law as
$(\bar{X})$. In other words, we will show that the law of the solution
$X_t(r)$, noted $m(t,r)(dy)$, is measurable with respect to
$\mathcal{B}(\Gamma)$, and that the mean-field equation can be
expressed as the integro-differential McKean--Vlasov equation
\begin{eqnarray*}
d \bar{X}_t(r)&=&f\bigl(r,t,\bar{X}_t(r)\bigr) \,dt
\\
&&{} + \int_{\Gamma} \int_E b\bigl(r,r',
\bar{X}_t(r),y \bigr) m\bigl(t-\tau\bigl(r,r'
\bigr),r'\bigr) (dy) \,d\lambda\bigl(r'\bigr) \,dt
\\
&&{} +
\sigma(r) \,dW_t(r),
\end{eqnarray*}
which will also be written, denoting $\mathcal{E}_{r'}$ is the
expectation with respect to the distribution of the population
locations over $\Gamma$ with distribution $\lambda(\cdot)$,
\begin{eqnarray*}
d \bar{X}_t(r)&=&f\bigl(r,t,\bar{X}_t(r)\bigr) \,dt +
\sigma(r) \,dW_t(r)
\\
&&{} + \mathcal{E}_{r'} \bigl[
\mathbb{E}_{\bar{Z}}\bigl[b\bigl(r,r',\bar {X}_t(r),
\bar{Z}_{t-\tau(r,r')}\bigl(r'\bigr)\bigr)\bigr] \bigr] \,dt.
\end{eqnarray*}
{\bf Update: } The spatially chaotic Brownian motion embodies the 
fact that neurons at different locations converge towards independent
processes. For each fixed location $r\in\Gamma$, the solution to the 
mean field equations has the same law $\check{m}(t,r)(dy)$ as the 
solution $\check{X}_{t}(r)$ to the equation:
\begin{eqnarray}\label{eq:MeasurableVersion}
d\check{X}_t(r)&=&f\bigl(r,t,\check{X}_t(r)\bigr) \,dt
+
\sigma(r) \,d\check{W}_t\\
\nonumber &&{} + \int_{\Gamma} \int_E b\bigl(r,r',
\check{X}_t(r),y \bigr) \check{m}\bigl(t-\tau\bigl(r,r'
\bigr),r'\bigr) (dy) \,d\lambda\bigl(r'\bigr) \,dt
\end{eqnarray}
where $\check{W}_{t}$ is a standard Brownian motion. This representation 
has the advantage to provide a version of the solution which is measurable
with respect to the space parameter $r$ (as soon as the initial condition is),
which thus simplifies mathematical developments. 

Let us eventually give the Fokker--Planck equation on the possible
density $p(t,r,y)$ of $m(r,t)$ with respect to Lebesgue's measure
%
\begin{eqnarray}\label{eqFP}
&& \partial_t p(t,r,x)\nonumber
\\
&&\qquad =-\nabla_x \biggl\{\biggl(f(r,t,x)
\nonumber\\[-8pt]\\[-8pt]
&&\hspace*{65pt}{} + \int_{\Gamma} \int_E b
\bigl(r,r',x,y \bigr) p\bigl(t-\tau\bigl(r,r'
\bigr),r',y\bigr) \,d\lambda\bigl(r'\bigr) \biggr)p(t,x)
\biggr\}\nonumber
\\
&&\quad\qquad{}+ \frac{1} 2 \Delta_x \bigl[ \bigl|\sigma(r)
\bigr|^2 p(t,x) \bigr].
\nonumber
\end{eqnarray}

The mean-field equations~(\ref{eqMFESpace}) are of a new type: they
resemble McKean--Vlasov equations but involve delays, spatially chaotic
Brownian motions and an ``integral over spatial locations.'' This is
hence a very unusual stochastic equation we need to thoroughly study in
order to ensure that these make sense and are well-posed. The existence
and uniqueness of solutions to these equations is addressed in
Section~\ref{secExistenceUniquenessSpace}, and the proof of the
propagation of chaos and convergence of the network equations toward
the solutions of that equations is addressed in
Section~\ref{secPropaChaSpace}.


\section{Analysis of the mean-field equation}\label{secExistenceUniquenessSpace} The mean-field
equation~(\ref{eqMFESpace}) involves two unusual terms: a stochastic
integral involving spatially chaotic Brownian motions and an integrated
McKean--Vlasov mean-field term.

Let us start by discussing properties of stochastic integrals with
respect to a spatially chaotic Brownian. Considering $\Delta_t(r)$ a
$\mathcal{F}_t$-progressively measurable process indexed by
$r\in\Gamma$ such that for any $r\in\Gamma$ we have
\begin{equation}
\label{eqBoundedDelta} \int_0^t \mathbb{E} \bigl[\bigl|
\Delta_s(r)\bigr|^2 \bigr] \,ds<\infty.
\end{equation}
It is trivial to see that for any $r\in\Gamma$, the process
$N_t(r):=\int_0^t \Delta_s(r)\,dW_s(r)$ is a well defined, square
integrable martingale with quadratic variation $\int_0^t |
\Delta_s(r)|^2 \,ds$.

The possible solutions $(\bar{X}_t(r))_{t,r}$ of the mean-field
equation have a law belonging to the set of probability measures on the
continuous functions of $[-\tau,T]$ with values in the set of mappings
of $\Gamma$ in $E$. It is important to note at this point that
similar to the spatially chaotic Brownian motion, the solutions of
the mean field equations are not measurable in $(\Gamma,\mathcal
{B}(\Gamma))$
since the solution considered at different space locations $r$ and $r'$
in $\Gamma$, namely $X_t(r)$ and $X_t(r')$, are independent.

Though trajectories of spatially chaotic processes are nonmeasurable,
their probability distribution, defining a set of measures parametrized
by $r\in\Gamma$, might be measurable. This is a necessary property to
make sense of the mean-field equations. Handling this subtlety
necessitates that we thoroughly define the space in which we are
working and where the mean-field equations are well defined. We define
$\mathcal{Z}$ the set of random variables whose law is measurable with
respect to $\mathcal{B}(\Gamma)$ [the random variable itself is not
assumed measurable with respect to $\mathcal{B}(\Gamma)$]. More
precisely, $\mathcal{Z}$ correspond to random variables whose law are
given by Markov kernels from $(\Gamma,\mathcal{B}(\Gamma))$ to\vadjust{\goodbreak}
$(\Omega,\mathcal{F})$, that is, mappings $p$ that associate to each
point $r\in\Gamma$ a probability measure $p(r)$ on
$(\Omega,\mathcal{F})$ such that for every measurable set
$A\in\mathcal{F}$, the map $r\mapsto p(r)(A)$ is measurable with
respect to $(\Gamma,\mathcal{B}(\Gamma))$. For a random variable
$(Z(r))_{r\in\Gamma}$ in $\mathcal{Z}$ with measurable law $p(r,dx)$,
we define with a slight abuse of notations the
$\mathbb{L}^2_{\lambda}(\Gamma)$ norm on $\Gamma$ by defining, for
$\hat{r}$ a $(\Omega',\mathcal{F}',\mathbb{P}')$ random variable with
law $\lambda$,
\begin{equation}
\label{eqNorm} \llVert Z \rrVert _{\mathbb{L}^2_{\lambda}(\Gamma
)}^2=\mathcal{E}_{\hat{r}} \bigl[\mathbb{E}
\bigl[\bigl| Z(\hat{r})\bigr| ^2\bigr] \bigr]
= \int_{\Gamma}\int_{E} x^2
p(r,dx)\,d\lambda(r),
\end{equation}
where $\mathcal{E}$ denotes the expectation on $\Omega'$. This clearly
defines a norm on random variables indexed by $r\in\Gamma$, when
identifying processes that are $\lambda$-a.e. $\mathbb{P}$-a.s. equal.
We denote $\mathbb{L}^2_{\lambda}(\Gamma)$ the set of random variables
in $\mathcal{Z}$ such that $\llVert  Z \rrVert
_{\mathbb{L}^2_{\lambda}(\Gamma)}<\infty$.

\begin{example*}
(i) The spatially chaotic Brownian motion at
fixed time $t$ has, for all $r\in\Gamma$, the
law of a standard Brownian motion. This law,
independent of $r\in\Gamma$, is hence
measurable with respect to $(\Gamma,\mathcal{B}(\Gamma))
$. Moreover, it belongs to $\mathbb{L}^2_{\lambda}(\Gamma)$ and has a~norm
equal to $t$.

(ii) Another example is given by the variable $Z_T(r)=\int_0^T
    \Delta_s(r)\,dW_s(r)$ where $\Delta$ is a function of
    $\mathbb{R}^+\times\Gamma$ measurable with respect to the
    $\sigma$-algebra
    $\mathcal{B}(\mathbb{R}^+)\otimes\mathcal{B}(\Gamma)$ and
    satisfying the condition $\int_0^T\int_{\Gamma}|
    \Delta_s(r)|^2 \,d\lambda(r) \,ds$. The thus defined
    variable is not measurable with respect to
    $\mathcal{B}(\Gamma)$, but belongs to $\mathcal{Z}$ since this
    variable is a centered Gaussian process with measurable
    variance $\int_0^T |\Delta_s(r)|^2\,ds$, hence the law
    of $Z(r)$ is $\mathcal{B}(\Gamma)$-measurable. Eventually,
    $Z\in\mathbb{L}^2_{\lambda}(\Gamma)$ with $\llVert  Z \rrVert
    _{\mathbb{L}^2_{\lambda}(\Gamma)}=\int_0^T
    \int_{\Gamma}|\Delta_s(r)|^2\,d\lambda(r)\,ds$.
\end{example*}
It is clear that for any random variable $Z(r)$ measurable in law, one can find
\footnote{A natural choice is 
to consider the variable $F^{{-1}}(r,U)$ where $F(r,\cdot)$ is the repartition 
function of $Z(r)$ and $U$ is a uniform random variable.} 
a process $\check{Z}_{t}(r)$ that is measurable with respect to $r\in\Gamma$. 

We extend this norm to processes with values in $\mathbb{L}^2_{\lambda
}(\Gamma)$. For $(Z_t(r))_{t\in[u,v]}$ a~stochastic process with
continuous paths indexed by $r\in\Gamma$ such that the law of $Z_t(r)$
is measurable with respect to $\mathcal{B}(\Gamma)$, we say that it 
belongs to $\mathcal{M}:=\mathcal{M}^2([u,v],\mathbb{L}^2_{\lambda}(\Gamma))$ 
if there exists a process $\check{Z}_t(r)$ measurable with respect to 
the product space $(\Omega \times \Gamma)$ such that for all fixed $r\in\Gamma$,
$Z_{t}(r)$ and $\check{Z}_{t}(r)$ have the same law, and moreover
\[
\| Z \|_{\mathcal{M}}:= \mathcal{E}_{\hat{r}} \Bigl(\mathbb {E}
\Bigl[\sup_{s\in[u,v]}\bigl|{\check{Z}_s(\hat{r})}
\bigr|^2 \Bigr] \Bigr)<\infty
\]
and this quantity defines a norm on $\mathcal{M}^2([u,v],\mathbb
{L}^2_{\lambda}(\Gamma))$ where are
identified the processes that are $\lambda$-a.e. and $\mathbb
{P}$-a.s. equal for all times. We will sometimes denote this norm
$\mathcal{E}_{\hat{r}} \Bigl(\mathbb {E}
\Bigl[\sup_{s\in[u,v]}\bigl|{Z_s(\hat{r})}
\bigr|^2 \Bigr] \Bigr)$ with a slight abuse of notations.

\begin{example*}
(i) The spatially chaotic Brownian motion on $[0,T]$ belongs to
    $\mathcal{M}^2([0,T], \mathbb{L}^2_{\lambda }(\Gamma))$ and has
    a norm equal to $T$ thanks to the classical property that the
    supremum of the Brownian motion has the law of the absolute
    value of the Brownian motion.

(ii) The process $Z_t(r)=\int_0^t \Delta_s(r)\,dW_s(r)$ introduced
    above belongs to
    $\mathcal{M}:=\mathcal{M}^2([0,T],\mathbb{L}^2_{\lambda}(\Gamma))$, since
    for all $r\in\Gamma$, $(Z_{t}(r))_{{t\geq 0}}$ has the same law as 
    $\check{Z}_{t}(r)=\int_0^t \Delta_s(r)\,dB_s$ for $B_{s}$ a standard 
    Brownian motion, process which is clearly measurable with respect to $r$.
    Moreover, thanks to Burkholder--Davis--Gundy inequality, has a~norm
    $\|{Z}\|_{\mathcal{M}}\leq4 \int_0^T \int_{\Gamma}|
    \Delta_s(r)|^2\,d\lambda(r)\,ds$.
\end{example*}

We will also work with processes $(X_{t}(r))_{t\geq 0, r\in\Gamma}$ that are measurable on the product space $(\Omega \times \Gamma)$. We denote $\mathcal{M}'=\mathcal{M}^{2}([-\tau,T]\times \Gamma,E)$ the space of square integrable processes, i.e. such that:
\[
\| X \|_{\mathcal{M'}}:= \int_{\Gamma} \mathbb {E}
\Bigl[\sup_{s\in[u,v]}\bigl|{X_s(r)}
\bigr|^2 \Bigr] d\lambda(r)<\infty
\]
Now that these norms are introduced, we are in position to show the
well-posedness of the mean-field equations:

\begin{theorem}\label{thmExistenceUniquenessSpace}
For any $(\zeta^0_t(r), t\in[-\tau,0], r\in\Gamma) \in
\mathcal{M}^2([-\tau,0],\mathbb{L}^2_{\lambda}(\Gamma))$\break  a~square-integrable process, the mean-field equation (\ref{eqMFESpace})
with initial condition $\zeta^0$ has a unique solution on
$[-\tau,T]$ for any $T>0$.
\end{theorem}

{\bf Update:} It is in the proof of this theorem that we make use of the 
$\Gamma$-measurable process $(\check{X}_{t}(r))_{t\geq 0}$ solution of (\ref{eq:MeasurableVersion})
and having the same law as $(\bar{X}_{t}(r))_{t\geq 0}$, to justify the fact that
all estimates used in the proof are well-defined.

\begin{pf}
As always for these types of properties, we reduce the problem to the
existence and uniqueness of a fixed point of a map $\Phi$ acting on
stochastic processes $X$ in $\mathcal{M}$
\begin{eqnarray*}
{\Phi}(X)_t(r) &=& \cases{\displaystyle{\zeta}^0_0(r) + \int
_0^t f\bigl(r,s,X_s(r)\bigr) \,ds+
\int_0^t \sigma(r) \,dW_s(r)
\vspace*{3pt}\cr
\qquad{} + \displaystyle\int_0^t \int_{\Gamma}
\mathbb{E}_{Z} \bigl[ b\bigl(r,r',X_s(r),
Z_{s-\tau(r,r')} \bigl(r'\bigr) \bigr)\bigr]\,d\lambda
\bigl(r'\bigr) \,ds,
\cr
\hspace*{79.5pt} t>0,
\vspace*{2pt}\cr
\zeta^0_t(r),\hspace*{50pt} t\in[-\tau, 0],
\vspace*{2pt}\cr
(Z_t)\stackrel{\mathcal{L}}{=}(X_t),\qquad\mbox{independent of }(X_t)\mbox{ and }\bigl(W_t(\cdot)\bigr).}
\end{eqnarray*}
Similarly, for $X\in \mathcal{M}'$, we define the map:
\begin{eqnarray*}
\check{\Phi}(X)_t(r) &=& \cases{\check{\zeta}^0_0(r) + \int
_0^t f\bigl(r,s,X_s(r)\bigr) \,ds+
\int_0^t \sigma(r) \,dW_s
\vspace*{3pt}\cr
\qquad{} + \displaystyle\int_0^t \int_{\Gamma}
\mathbb{E}_{Z} \bigl[ b\bigl(r,r',X_s(r),
Z_{s-\tau(r,r')} \bigl(r'\bigr) \bigr)\bigr]\,d\lambda
\bigl(r'\bigr) \,ds,
\cr
\hspace*{79.5pt} t>0,
\vspace*{2pt}\cr
\check{\zeta}^0_t(r),\hspace*{50pt} t\in[-\tau, 0],
\vspace*{2pt}\cr
(Z_t)\stackrel{\mathcal{L}}{=}(X_t),\qquad\mbox{independent of }(X_t)\mbox{ and }\bigl(W_t(\cdot)\bigr).}
\end{eqnarray*}
%
We can show that for $X\in \mathcal{M}$, $\Phi(X)$ is in $\mathcal{M}$. It is clear that $\Phi(X)$ is a spatially chaotic process. If $\check{X}$ is a $(\Omega\times \Gamma)$-measurable process such that for all fixed $r$, $(X_{t}(r))_{t\geq -\tau}$ has the same law as $(\check{X}_{t}(r))_{t\geq -\tau}$, then $(\Phi(X)_{t}(r))_{{t\geq 0}}$ has the same law as $(\check{\Phi}(\check{X})_{t}(r))_{{t\geq 0}}$. Moreover,  $(\check{\Phi}(\check{X})_{t}(r))$ is clearly $\Omega\times\Gamma$-measurable, and square integrable, since:
%
\begin{eqnarray}
\label{eqBoundPhi}
\bigl\llVert \Phi(X) \bigr\rrVert _{\mathcal{M}} &\leq& 4 \biggl(
\bigl\llVert \zeta_0^0 \bigr\rrVert _{\mathbb{L}^2_{\lambda
}(\Gamma)} + T C
\int_0^T \bigl(1+\llVert X
\rrVert _{\mathcal{M}}\bigr)\nonumber
\nonumber\\[-8pt]\\[-8pt]
&&\hspace*{11pt}{} +4 \int_0^T \int
_{\Gamma}\bigl|\sigma (r)\bigr|^2\,d\lambda(r) \,ds
+ T\tilde{K}\int_0^T
\bigl(1+\llVert X \rrVert _{\mathcal{M}}^2\bigr) \,ds
\biggr),\hspace*{-22pt}\nonumber
\end{eqnarray}
which is finite for any $X\in\mathcal{M}$.

We may hence iterate the map $\Phi$. We fix $X^{0}$ a process in
$\mathcal{M}$ and build the sequence $X^k$ by induction through the
recursion relationship $X^{k+1}=\Phi(X^k)$. We aim to show that these
processes constitute\vadjust{\goodbreak} a Cauchy sequence in $\mathcal{M}$
To this purpose, we
introduce the norm of the process up to time $t$
\[
\| X \|_{\mathcal{M}_t}^2:= \mathcal
{E}_r \Bigl[ \mathbb{E}\Bigl(\sup_{s\in[-\tau,t]}\bigl|{X_s}(r)
\bigr|^2 \Bigr) \Bigr]= \int_{\Gamma}\mathbb{E}\biggl[
\sup_{s\in[-\tau,t]}\bigl| {X_s}(r)\bigr|^2 \biggr]\,d
\lambda(r).
\]
We now introduce a sequence of processes $(Z^k)$ independent of the
collection of processes $(X^k)$ and having the same law, built
recursively as follows:
\begin{itemize}
\item $Z^0$ is independent of $X^0$ and has the same law as $X^0$;

\item for $k\geq1$, $Z^{k}$ is independent of the sequence
of processes $(X^{0},\ldots,X^k)$ and is such that the
collection of processes $(Z^{0},\ldots,Z^{k})$ has the
same joint law as $(X^{0},\ldots,X^k)$, that is, $Z^k$
is chosen such as its conditional law given
$(Z^0,\ldots,Z^{k-1})$ is the same as that of $X^k$
given $(X^0,\ldots,X^{k-1})$.
\end{itemize}
We study the norm $\| X^{k+1}-X^k \|_{\mathcal{M}_T}$. We
decompose this
difference into the sum of three elementary terms as follows, for
$t\in[0,T]$ and $r\in\Gamma$:
\begin{eqnarray*}
X^{k+1}_t(r)-X^{k}_t(r) &=& \int
_0^t \bigl\{ \bigl(f\bigl(r,s,X_s^{k}(r)
\bigr)-f\bigl(r,s,X_s^{k-1}(r)\bigr)\bigr) \bigr\} \,ds
\\
&&{}+ \int_0^t \int_{\Gamma}
\bigl\{ \bigl( \mathbb{E}_{Z} \bigl[ b\bigl(r,r',X_s^{k}(r),
Z^{k}_{s-\tau(r,r')} \bigl(r'\bigr) \bigr)\bigr]
\\
&&\hspace*{45pt}{} -\mathbb{E}_{Z} \bigl[ b\bigl(r,r',X^{k-1}_s(r),
Z^{k}_{s-\tau(r,r')} \bigl(r'\bigr) \bigr)\bigr] \bigr)
\bigr\}\,d\lambda\bigl(r'\bigr) \,ds
\\
&&{}+ \int_0^t \int_{\Gamma}
\bigl\{ \bigl( \mathbb{E}_{Z} \bigl[ b\bigl(r,r',X_s^{k-1}(r),
Z^{k}_{s-\tau(r,r')} \bigl(r'\bigr) \bigr)\bigr]
\\
&&\hspace*{45pt}{}-\mathbb{E}_{Z} \bigl[ b\bigl(r,r',X^{k-1}_s(r),
Z^{k-1}_{s-\tau(r,r')} \bigl(r'\bigr) \bigr)\bigr] \bigr)
\bigr\}\,d\lambda\bigl(r'\bigr) \,ds
\\
& =:& A_t(r) + B_t(r)+ C_t(r).
\end{eqnarray*}
We hence obviously have
\[
M^k_t:=\bigl\| X^{k+1}-X^{k}
\bigr\|_{\mathcal{M}_t}^2\leq3 \bigl(\| A \|_{\mathcal{M}_t}^2
+ \| B \|_{\mathcal{M}_t}^2+ \| C \|_{\mathcal{M}_t}^2
\bigr).
\]
We treat each term separately, and denote $\check{X}^{k}=\check{\Phi}^{k}(\check{X}^{0})$ the $(\Omega\times \Gamma)$ 
measurable version of $X^{k}$. We have
\begin{eqnarray*}
\llVert A \rrVert _{\mathcal{M}_{t}}^2 &=& \mathbb{E} \biggl[\int
_{\Gamma} \biggl( \sup_{s\in[0,t]}\biggl|\int
_0^s f\bigl(r,u,\check{X}_u^{k}(r)
\bigr)
\\
&&\hspace*{56pt}{} -f\bigl(r,u,\check{X}_u^{k-1}(r)\bigr) \,du \biggr|^2
\biggr)\,d\lambda(r) \biggr]
\\
\mbox{(Cauchy--Schwarz)} &\leq& T K_f^2 \mathbb{E}
\biggl[\int_{\Gamma
} \biggl(\int_0^t
\bigl| \check{X}_s^{k}(r)-\check{X}_s^{k-1}(r)
\bigr|^2\,ds \biggr)\,d\lambda(r) \biggr]
\\
&\leq& T K_f^2 \int_0^t
\bigl\llVert {X}^{k}-{X}^{k-1} \bigr\rrVert _{\mathcal{M}_{s}}^2
\,ds
\end{eqnarray*}
which directly implies $\llVert  A \rrVert _{\mathcal {M}_{t}}^2\leq T
K_f^2 \int_0^t M_s^{k-1} \,ds$. The terms $B_t$ and $C_t$ can be
controlled using the same techniques. Let us, for instance, treat the
case of $C_t$. We have
\begin{eqnarray*}
\llVert C \rrVert _{\mathcal{M}_{t}}^2 & =& \mathbb{E} \biggl[\int
_{\Gamma}\sup_{s\in[0,t]} \biggl|\int
_{\Gamma}\int_0^s \bigl(
\mathbb{E}_{Z} \bigl[ b\bigl(r,r',\check{X}_u^{k-1}(r),
Z_{u-\tau(r,r')}^{k} \bigl(r'\bigr) \bigr)\bigr]
\\
&&\hspace*{81pt}{}- \mathbb{E}_{Z} \bigl[ b\bigl(r,r',\check{X}_u^{k-1}(r),
\\
&&\hspace*{126pt}{}
Z_{u-\tau (r,r')}^{k-1} \bigl(r'\bigr) \bigr)\bigr] \bigr)
\,du \,d\lambda\bigl(r'\bigr) \biggr|^2 \,d\lambda(r) \biggr]
\\
\mbox{(CS)} & \leq& t \int_{\Gamma^2} \int_0^t
\mathbb{E} \bigl[ \mathbb {E}_{Z} \bigl[ \bigl| b\bigl(r,r',\check{X}_u^{k-1}(r),
Z_{u-\tau(r,r')}^{k} \bigl(r'\bigr) \bigr)
\\
&&\hspace*{67pt}{} - b \bigl(r,r',\check{X}_u^{k-1}(r), Z_{u-\tau(r,r')}^{k-1}
\bigl(r'\bigr) \bigr)\bigr| ^2\bigr] \bigr) \,du \bigr]\,d
\lambda(r)\,d\lambda\bigl(r'\bigr)
\\
\mbox{\hyperref[AssumpLocLipschbSpace]{(H2)}}&\leq& t L^2 \int_{\Gamma^2} \int
_0^t \mathbb{E} \bigl[\bigl| \check{X}_s^{k}(r)-\check{X}_s^{k-1}(r)
\bigr|^2 \bigr] \,ds \,d\lambda(r)\,d\lambda\bigl(r'\bigr)
\\
&\leq& t L^2 \int_0^t \mathbb{E}
\bigl[\bigl\llVert \check{X}_s^{k}(r)-\check{X}_s^{k-1}(r)
\bigr\rrVert _{\mathbb{L}^2_{\lambda}(\Gamma
)}^2 \bigr] \,ds = t L^2 \int
_0^t M_s^{k-1} \,ds.
\end{eqnarray*}
The term $B_t$ is treated exactly in the same manner and yields the inequality
\[
\llVert B \rrVert _{\mathcal{M}_{t}}^2\leq t L^2 \int
_0^t M_s^{k-1} \,ds.
\]
All together we obtain, using the fact that for all $k>1$,
$t\in[-\tau,0]$ and $r\in\Gamma$ we have $X^{k+1}_t(r) = X^{k}_t(r) =
\zeta^0_t(r)$
\begin{equation}
\label{eqBoundaryMRec} M_t^k\leq K' \int
_{0}^t M_s^{k-1} \,ds
\end{equation}
with $K'=3 T (K_f^2+2 L^2)$, readily implying
\begin{equation}
\label{eqCauchySeq} M^k_t \leq\bigl(K'
\bigr)^k \int_0^t\int
_0^{s_1}\cdots\int_0^{s_{k-1}}M^0_{s_k}
\,ds_1\cdots ds_k \leq\frac{(K')^k t^k}{k!}
M^0_T
\end{equation}
and $M^0_t$ is finite since we assumed $X\in\mathcal{M}_T$ and showed
that $\Phi(X)\in\mathcal{M}_T$. Routine methods starting from
inequality~(\ref{eqCauchySeq}) using the Benaym\'e--Chebychev
inequality and the Borel--Cantelli lemma allow us to prove existence
and uniqueness of a fixed point $X^{*}\in\mathcal{M}$ for $\Phi$ (see, e.g.,
\cite{revuz-yor99}, pages 376--377), and that this fixed point is
adapted and almost surely continuous.
The same convergence holds for the sequence $\check{\Phi}$, and the fixed point of
$\Phi$ has, for all fixed location, the same law as $\check{X}^{*}$ the unique fixed 
point of $\hat{\Phi}$. This process being limit of a sequence of 
$(\Gamma,\mathcal{B}(\Gamma))$-measurable processes, 
it is measurable as well, and since it satisfies $\check{\Phi}(\check{X}^{*})=\check{X}^{*}$. Using inequality~(\ref{eqBoundPhi}) and the fact that
${X}^{*}=\Phi({X}^{*})$ we have
\[
\llVert {X}^{*} \rrVert _{\mathcal{M}_{t}} \leq4 \biggl(\bigl\llVert
\zeta_0^0 \bigr\rrVert _{\mathbb{L}^2_{\lambda
}(\Gamma)} + T \bigl(C+
\tilde{K}+4\llVert \sigma\rrVert _{\mathbb{L}^2_{\lambda}(\Gamma)}\bigr) + T(C+\tilde{K}) \int
_0^{t} \llVert {X}^{*} \rrVert _{\mathcal{M}_{s}} \,ds
\biggr)
\]
ensuring by Gronwall's lemma that the solution has a finite norm in
$\mathcal{M}_T$. Proving uniqueness of the solution using
equation~(\ref{eqBoundaryMRec}) is then folklore.

\end{pf}

Now that we proved existence and uniqueness of solutions for the
mean-field equations, we now turn to showing that network equations
indeed converge in law toward this solution and that the propagation of
chaos occurs.

\section{Limit in law and propagation of chaos}\label{secPropaChaSpace}
We are now in a position to prove the main result of the manuscript,
namely the convergence in law of the solutions of the network
equations~(\ref{eqNetworkSpace})--(\ref{eqMFESpace})
and the fact that the propagation of chaos property occurs. To this
end, we consider that the network equations have chaotic initial
conditions. In detail, let $(\zeta^0_t(r)) \in\mathcal{M}^2([-\tau,0],
\mathbb{L}^2_{\lambda}(\Gamma))$ a spatially chaotic stochastic
process, that is, a stochastic process such that for any $r\neq r'$,
the process $(\zeta^0_t(r))$ is independent of $(\zeta^0_t(r'))$. We
consider that the initial condition of different neurons in the network
are independent and the initial condition
$(\zeta^{i}_t)\in\mathcal{M}^2([-\tau,0], \mathbb
{L}^2_{\lambda}(\Gamma))$ for neuron $i$ in population $\alpha$, is
equal to $(\zeta^{0}_t(r_{\alpha }))\in
\mathcal{M}^2(\mathcal{C}_{\tau})$.

The classical coupling argument cannot be directly applied here.
Indeed, the usual argument is based on the fact that we are able to
define the solution of the mean-field equation through the use of the
\emph{same} Brownian motion and with the \emph{same} initial condition
as one of the neurons (or particles). This is no more~the case because
individual neurons are governed by finite-dimensional Brownian motions
and the mean-field equation by a spatially chaotic Brownian motion.
Notwithstanding, an argument based on a slightly more subtle couplings
holds. In detail, let us consider neuron $i\in\mathbb{N}$ of the
network, in
population $\alpha$ at location $r_{\alpha} \in\Gamma$. Denote by
$(\tilde{W}^{i}_t)$ the Brownian motion governing the evolution of
neuron $i$ in the network and $\zeta^i\in\mathcal{M}(\mathcal
{C}_{\tau})$ the initial condition
of the network. We aim at defining a spatially chaotic Brownian motion
$W^i_t(r)$ on $\mathbb{R}^{m\times d}$ such that the standard Brownian motion
$(W^i_t(r_{\alpha}))$ is equal to $(\tilde{W}^i_t)$, and proceed as
follows. Let $(W_t(r))_{t\in[0,T], r\in\Gamma}$ be a ${m\times
d}$-dimensional spatially chaotic Brownian motions independent of the
processes $(\tilde{W}_t^{j})$. The processes
\[
\cases{ \bigl(W^i_t(r)\bigr) = \bigl(W_t(r)
\bigr), &\quad $r\neq r_{\alpha}$,
\vspace*{2pt}\cr
\bigl(W^i_t(r_{\alpha})
\bigr) = \bigl(\tilde{W}^i_t\bigr)}
\]
are clearly spatially chaotic Brownian motions and will be used to
construct a~particular solution of the mean-field equations. In order
to completely define a solution of\vadjust{\goodbreak} the mean-field equations, we need to
specify an initial condition, and aim at coupling it to the initial
condition of neuron $i$. To this end, we define a spatially chaotic
process $(\tilde{\zeta}^0_t(r)) \in\mathcal{M}^2([-\tau,0],
\mathbb{L}^2_{\lambda}(\Gamma))$ equal in law
to $(\zeta^0_t(r))$ and independent of $\zeta^i_t$, and define a
coupled process $({\zeta}^{i,0}_t(r)) \in\mathcal{M}^2([-\tau,0],
\mathbb{L}^2_{\lambda}(\Gamma))$~as
\[
\cases{ \zeta^{i,0}_t(r) = \tilde{\zeta}^0_t(r),
&\quad$r\neq r_{\alpha
}$,
\vspace*{2pt}\cr
\zeta^{i,0}_t(r_{\alpha})
= \zeta^i_t.}
\]
Here again, it is clear that this process is spatially chaotic, that
is, that for any $r\neq r'$, the processes $\zeta^{i,0}_t(r)$ and
$\zeta^{i,0}_t(r')$ are independent, and that $\zeta^{i,0}_t(r)$ has
the law of $\zeta^0_t(r)$.

Now that these processes have been constructed, we are in a position to
define the process $(\bar{X}^i_t)$ as the unique solution of mean-field
equation~(\ref{eqMFESpace}), driven by the spatially chaotic Brownian
motion $(W_t^i(r))$ and with the spatially chaotic initial condition
$(\zeta^{i,0}_t(r))$
\[
\cases{\displaystyle d\bar{X}^i_t(r) = f\bigl(r,t,
\bar{X}^i_t(r)\bigr) \,dt + \int_{\Gamma}
\mathbb{E}_Z\bigl[b\bigl(r,r',\bar{X}^i_t(r),Z_{t-\tau(r,r')}
\bigl(r'\bigr)\bigr)\bigr]\,d\lambda \bigl(r'\bigr) \,dt
\cr
\hspace*{43pt}{} + \displaystyle\sigma(r) \,dW^i_t(r),\qquad t\geq0,
\vspace*{2pt}\cr
\displaystyle\bar{X}^i_t(r) = \zeta^{i,0}_t (r),\hspace*{65.5pt}  t\in[-\tau, 0],
\vspace*{2pt}\cr
\displaystyle (Z_t)\stackrel{\mathcal{L}} {=}
\bigl(\bar{X}^i_t\bigr) \in\mathcal{M},\hspace*{57pt} \mbox{independent of }\bigl(\bar{X}^i_t\bigr),
\bigl(W^{i}_t(\cdot)\bigr)\mbox{ and }
\bigl(B^{i}_t(\cdot,\cdot)\bigr).}
\]
The same procedure applied to all $j\in\mathbb{N}$ allows us to build a
collection of independent stochastic processes
$(\bar{X}^j_t(r))_{j=1,\ldots, N} \in\mathcal{M}^2([-\tau,T],\mathbb
{L}^2_{\lambda}(\Gamma))$ such that all neurons $j$ in population
$\alpha$ have the same law as $(\bar{X}(r_{\alpha}))$. Let us denote by
$m(t,r)$ the probability distribution of $\bar{X}_t(r)$ solution of the
mean-field equation~(\ref{eqMFESpace}). As previously, the process
$(Z_t(r))$ generically denotes a process belonging to
$\mathcal{M}^2([-\tau,T],\mathbb {L}^2_{\lambda}(\Gamma))$ and
distributed as $m$.

Let us fix $l\in\mathbb{N}^*$ and $(i_1,\ldots,i_l)$, a collection of
neuron indexes, respectively, belonging to populations located at
$(r_1, \ldots, r_k)$ (possibly identical). We now prove the almost sure
convergence of a collection of processes $({X}^{i_k,N}_t, k=1,\ldots,
l)$ toward $(\bar{X}^{i_k}_t(r_k), k=1,\ldots,l)$, implying its
convergence of the law toward the chaotic distribution $m {(t,r_1)}
\otimes\cdots\otimes m(t, r_k)$ as $N$ goes to infinity. We start by
proving this property for $l=1$ before extending that result to $l> 1$.

\begin{theorem}\label{thmPropagationChaosSpace}
Let $i\in\mathbb{N}$ a fixed neuron in population $\alpha$. Under
assumptions~\textup{\hyperref[AssumpLocLipschSpace]{(H1)}--\hyperref[AssumpSpaceContinuity]{(H5)}}
and the neural field assumption~(\ref{eqPopulationEstimate}), for
almost all realizations of the population locations $(r_{\alpha},
\alpha\in\mathbb{N})$, the process $(X^{i,N}_t, t\leq T)$ solution of
the network equations (\ref{eqNetworkSpace}) converges in law toward
the process $(\bar{X}_t(r_{\alpha}), t\leq T)$ solution of the
mean-field equations (\ref{eqMFESpace}) with initial condition
$(\zeta^0_t(r))$, and moreover, the speed of convergence is given by
\begin{equation}
\label{eqPropchaosSpace} \mathcal{E} \Bigl( \mathbb{E} \Bigl[\sup_{-\tau\leq s\leq T}
\bigl| X^{i,N}_s - \bar{X}^i_s(r_{\alpha})
\bigr|^2 \Bigr] \Bigr) = O \biggl(\mathbbe (N) + \frac{1}{P(N)}
\biggr).\vadjust{\goodbreak}
\end{equation}
\end{theorem}

\begin{remark*}
We recall that $\mathcal{E}$ denotes the expectation on the
distribution of the space locations $(r_k)_{k=1,\ldots,P(N)}$ and
$\mathbbe(N)=\frac{1} {P(N)} \sum_{\gamma=1}^{P(N)}\frac{1} {
N_{\gamma}(N)}$.
\end{remark*}

\begin{pf}
The proof is based on evaluating the distance\break  $\mathbb{E}[\sup_{-\tau
\leq
s\leq T} | X^{i,N}_s - \bar{X}^i_s|^2 ]$, and breaking it into
a few elementary, easily controllable terms. A substantial difference
with usual mean-field proofs is that we need to prove a convergence in
the infinite-dimensional space $\mathbb{L}^2_{\lambda}(\Gamma)$, and
that the interaction term in
networks equations consists of a sum over a finite number of
populations, whereas the effective interaction term arising in the
mean-field equation is an integral over $\Gamma$.

Throughout the demonstration, we will generically denote by
$r_{\beta}\in\Gamma$ the location of population
$\beta\in\{1,\ldots,P(N)\}$. We use the following elementary
\mbox{decomposition} [each line of the righthand side corresponds to one term
of the decomposition, $A_t(N)-E_t(N)$]:
\begin{eqnarray*}
&& X^i_t-\bar{X}^i_t(r_{\alpha})
\\
&&\qquad = \int_0^t \bigl(f\bigl(r_{\alpha
},s,X^i_s \bigr)-f\bigl(r_{\alpha},s,\bar{X}^i_s(r_{\alpha})
\bigr)\bigr) \,ds
\\
&&\quad\qquad{} + \frac{1}{P(N)} \sum_{\gamma=1}^{P(N)}
\int_0^t \frac{1} {
N_{\gamma}} \sum
_{j=1}^{N_{\gamma}} \bigl(b\bigl(r_{\alpha},r_{\gamma
},X^i_s,X^j_{s-\tau(r_{\alpha},r_{\gamma})}
\bigr)
\\
&&\hspace*{147pt}{} -b\bigl(r_{\alpha
},r_{\gamma},\bar{X}^i_s(r_{\alpha}),X^j_{s-\tau(r_{\alpha
},r_{\gamma})}
\bigr)\bigr) \,ds
\\
&&\quad\qquad{}+ \frac{1}{P(N)} \sum_{\gamma=1}^{P(N)}
\int_0^t \frac{1} {
N_{\gamma}} \sum
_{j=1}^{N_{\gamma}} \bigl(b\bigl(r_{\alpha},r_{\gamma},
\bar {X}^i_s(r_{\alpha}),X^j_{s-\tau(r_{\alpha},r_{\gamma})}
\bigr)
\\
&&\hspace*{147pt}{} - b\bigl(r_{\alpha},r_{\gamma},\bar{X}^i_s(r_{\alpha}),
\bar{X}^j_{s-\tau
(r_{\alpha},r_{\gamma})}(r_{\gamma})\bigr)\bigr) \,ds
\\
&&\quad\qquad{}+ \frac{1}{P(N)} \sum_{\gamma=1}^{P(N)}
\int_0^t \Biggl(\frac{1} {N_{\gamma}} \sum
_{j=1}^{N_{\gamma}} b\bigl(r_{\alpha},r_{\gamma
},\bar{X}^i_s(r_{\alpha}),\bar{X}^j_{s-\tau(r_{\alpha},r_{\gamma
})}(r_{\gamma})
\bigr)
\\
&&\hspace*{118pt}{} -\mathbb{E}_Z\bigl[ b\bigl(r_{\alpha},r_{\gamma},
\bar {X}^i_s(r_{\alpha}),Z_{s-\tau(r_{\alpha},r_{\gamma})}(r_{\gamma
})
\bigr)\bigr] \Biggr) \,ds
\\
&&\quad\qquad{}+ \frac{1}{P(N)} \sum_{\gamma=1}^{P(N)}
\int_0^t \biggl( \mathbb{E}_Z
\bigl[b\bigl(r_{\alpha},r_{\gamma},\bar{X}^i_s(r_{\alpha
}),Z_{s-\tau(r_{\alpha},r_{\gamma})}(r_{\gamma})
\bigr)\bigr]
\\
&&\quad\qquad\hspace*{81pt}{} -\int_{\Gamma} \mathbb{E}_Z\bigl[b
\bigl(r_{\alpha},r',\bar{X}_s(r_{\alpha}),Z_{s-\tau
(r_{\alpha},r')}
\bigl(r'\bigr)\bigr)\bigr]\,d\lambda\bigl(r'\bigr)
\biggr) \,ds
\\
&&\qquad =: A_t(N)+B_t(N)+C_t(N)+D_t(N)+E_t(N).\nonumber
\end{eqnarray*}
Due to the exchangeability of neurons belonging to the same population,
the probability distribution of these terms does not depend on the
particular neuron $i$ considered, but only on the population it belongs
to. The terms $A_t(N)$, $B_t(N)$ and $C_t(N)$ involve the Lipschitz
continuity of the functions involved, the term $D_t(N)$ correspond to
averaging effects (mean-field limit) at single populations levels and
the term $E_t(N)$ corresponds to the continuous limit.
The terms $A_t(N)$ through $C_t(N)$ are treated using the Lipschitz
continuity of the functions involved. Using Cauchy--Schwarz (CS)
inequalities, we easily obtain
\begin{eqnarray*}
\mathbb{E}\Bigl[\sup_{0\leq s\leq t} \bigl| A_s(N)
\bigr|^2\Bigr] & \leq& K_f^2 T \int
_0^{t} \mathbb{E}\Bigl[\sup_{-\tau\leq u\leq s}
\bigl| X_u^{i,N}-\bar{X}_u^i(r_{\alpha})
\bigr|^2 \Bigr] \,ds,
\\
\mathbb{E}\Bigl[\sup_{0\leq s\leq t} \bigl| B_s(N)
\bigr|^2\Bigr] & \leq& T L^2 \int_0^{t}
\mathbb{E}\Bigl[\sup_{-\tau\leq u\leq s}\bigl| X^{i,N}_u-
\bar {X}^i_u(r_{\alpha}) \bigr|^2\Bigr]
\,ds,
\\
\mathbb{E}\Bigl[\sup_{0\leq s\leq t} \bigl| C_s(N)
\bigr|^2\Bigr] & \leq& T L^2 \int_0^{t}
\max_{j=1,\ldots, N} \mathbb{E}\Bigl[\sup_{-\tau\leq u \leq
s}\bigl|
X^{j,N}_u-\bar{X}^j_u(r_{p(j)})
\bigr|^2\Bigr] \,ds.
\end{eqnarray*}
Let us, for instance, treat the case of $B_t (N)$,
\begin{eqnarray*}
&& \mathbb{E}\Bigl[\sup_{0 \leq s\leq t}\bigl| B_s(N)
\bigr|^2\Bigr]
\\
&&\qquad = \frac
{1}{P(N)^2} \mathbb{E} \Biggl[\sup
_{0 \leq s\leq t} \Biggl|\sum_{\gamma=1}^{P(N)}
\int_0^s \frac{1} {N_{\gamma}} \sum
_{j=1}^{N_{\gamma}} \bigl(b\bigl(r_{\alpha},r_{\gamma
},X^{i,N}_u,X^{j,N}_{u-\tau(r_{\alpha},r_{\gamma})}
\bigr)
\\
&&\hspace*{178pt}{}-b\bigl(r_{\alpha},r_{\gamma},\bar{X}^{i}_u,X^{j,N}_{u-\tau
(r_{\alpha},r_{\gamma})}
\bigr) \bigr) \,du \Biggr|^2 \Biggr]
\\
&&\hspace*{-1.5pt} \mbox{(CS)}\leq\frac{T}{P(N)} \sum_{\gamma=1}^{P(N)}
\int_0^t \frac{1} {N_{\gamma}} \sum
_{j=1}^{N_{\gamma}} \mathbb{E} \bigl[\bigl\llvert b
\bigl(r_{\alpha},r_{\gamma},X^{i,N}_s,X^{j,N}_{s-\tau(r_{\alpha
},r_{\gamma})}
\bigr)
\\
&&\hspace*{145pt}{} -b\bigl(r_{\alpha},r_{\gamma},\bar {X}^{i}_s,X^{j,N}_{s-\tau(r_{\alpha},r_{\gamma})}
\bigr) \bigr\rrvert ^2 \bigr] \,ds
\\
&&\hspace*{-1pt} \mbox{\hyperref[AssumpLocLipschbSpace]{(H2)}} \leq T L^2 \int_0^t
\mathbb {E} \bigl[\bigl\llvert X^{i,N}_s-
\bar{X}^i_s\bigr\rrvert ^2 \bigr] \,ds
\\
&&\qquad \leq T L^2 \int_0^t \mathbb{E}
\Bigl[\sup_{-\tau\leq u\leq s}\bigl| X^{i,N}_u-
\bar{X}^i_u \bigr|^2 \Bigr] \,ds.
\end{eqnarray*}
The mean-field term $D_t(N)$ involves the difference between an
empirical mean of a function of processes and an expectation term, and
all have bounded second moment thanks to
Theorem~\ref{thmExistenceUniquenessSpace} and
assumption~\hyperref[AssumpbBoundSpace]{(H3)}. We have, using a (CS)
inequality,
\begin{eqnarray*}
&& \mathbb{E}\Bigl[\sup_{0\leq s\leq t}\bigl| D_s(N)
\bigr|^2\Bigr]
\\
&&\qquad\quad  \leq\frac{T }{P(N)} \sum
_{\gamma=1}^{P(N)} \int_0^t
\mathbb {E}\Biggl[ \Biggl|\frac{1} {N_{\gamma}} \sum_{j=1}^{N_{\gamma}}
b\bigl(r_{\alpha},r_{\gamma},\bar{X}^i_s,
\bar{X}^j_{s-\tau(r_{\alpha
},r_{\gamma})}\bigr)
\\
&&\hspace*{128pt}{} - \mathbb{E}_Z\bigl[b
\bigl(r_{\alpha},r_{\gamma},\bar {X}^i_s,Z_{s-\tau(r_{\alpha},r_{\gamma})}^{\gamma}
\bigr)\bigr] \Biggr| ^2\Biggr] \,ds
\end{eqnarray*}
and hence involves an expectation of the following type:
\begin{eqnarray*}
&& \mathbb{E}\Biggl[\Biggl|\frac{1} {N_{\gamma}} \sum_{j=1}^{N_{\gamma}}
\Theta\bigl(\bar{X}^i_s,\bar{X}^j_s
\bigr)-\mathbb{E}_Z\bigl[\Theta\bigl(\bar {X}^i_s,Z_s^{\gamma}
\bigr)\bigr]\Biggr|^2\Biggr]
\\
&&\qquad = \frac{1} {N_{\gamma}^2} \sum
_{k,l
=1}^{N_{\gamma}} \mathbb{E} \bigl[\bigl(\Theta\bigl(
\bar{X}^i_s,\bar {X}^j_s\bigr)-
\mathbb{E}_Z\bigl[\Theta\bigl(\bar{X}^i_s,Z_s^{\gamma}
\bigr)\bigr]\bigr)^T
\\
&&\hspace*{82pt}{}\cdot
\bigl(\Theta\bigl(\bar{X}^i_s,
\bar{X}^k_s\bigr)-\mathbb{E}_Z\bigl[\Theta
\bigl(\bar{X}^i_s,Z_s^{\gamma}\bigr)
\bigr]\bigr) \bigr],
\end{eqnarray*}
where $\Theta(x,y)=b(r_{\alpha},r_{\gamma},x,y)$. Routine methods allow
us to show that all the terms of the sum corresponding to indexes $j$
and $k$ such that the three conditions $j\neq i$, $k\neq i$ and $j\neq
k$ are satisfied are null. One simple way to show this property
consists of writing the expectations as integrals with respect to the
measure $m(t,r_{\alpha})$ and observing that all terms annihilate.
Therefore, there are no more than~$3 N_{\gamma}$ nonnull terms in the
sum (in the case $\alpha=\gamma$ there are just $N_{\gamma}$ nonnull
terms), and moreover, all of these terms are uniformly bounded. The
terms related to indexes $j =k \neq i $ satisfy the inequality
\begin{eqnarray*}
&& \mathbb{E} \bigl[\bigl\llvert \Theta\bigl(\bar{X}^i_s,
\bar {X}^j_s\bigr)-\mathbb{E}_Z\bigl[\Theta
\bigl(\bar{X}^i_s,Z_s^{\gamma}\bigr)
\bigr]\bigr\rrvert ^2 \bigr]
\\
&&\qquad \leq 2 \mathbb{E} \bigl[\bigl\llvert
\Theta\bigl(\bar {X}^i_s,\bar{X}^j_s
\bigr)\bigr\rrvert ^2 + \bigl\llvert \mathbb{E}_Z\bigl[
\Theta \bigl(\bar{X}^i_s,Z_s^{\gamma}
\bigr)\bigr]\bigr\rrvert ^2 \bigr]
\\
&&\qquad \leq 2 \Bigl\{\tilde{K}\bigl(1+\mathbb{E}\bigl[\bigl|\bar{X}^{i}_s
\bigr|^2\bigr]\bigr) + \mathbb{E} \Bigl[\Bigl\llvert
\mathbb{E}_Z\Bigl[\sqrt{\tilde{K}\bigl(1+\bigl|
\bar{X}^i_s\bigr|^2\bigr)}\Bigr] \Bigr\rrvert
^2 \Bigr] \Bigr\}
\\
&&\qquad \leq 4 \tilde{K} \bigl(1+C'(s)\bigr)
\end{eqnarray*}
with $C'(s)$ given by Theorem~\ref{thmExistenceUniquenessSpace}. The
terms related to the cases $j=i$ (or symmetrically $k=i$) are bounded
by the same constant, since we have for all $k$ such that $p(k)=\alpha
$, using the Cauchy--Schwarz inequality. We note $C= 4 \tilde{K}
(1+C'(T))$. We hence conclude that
\begin{eqnarray*}
\mathbb{E}\Bigl[\sup_{0\leq s\leq t}\bigl| D_s(N)
\bigr|^2\Bigr] &\leq& T^2 C \frac{1}{P(N)} \sum
_{\gamma=1}^{P(N)}\frac{3N_{\gamma
}-1}{N_{\gamma}^2}
\\
&\leq& 3 T^2 C
\frac{1}{P(N)} \sum_{\gamma
=1}^{P(N)}
\frac{1}{N_{\gamma}}= 3T^2 C \mathbbe (N).
\end{eqnarray*}

It hence only remains to control the term $E_t(N)$ corresponding to the
difference between an integral over the space $\Gamma$ weighted by the
density $d\lambda(r)$ and a sum, weighted by $1/P(N)$ of the same
integrand at $P(N)$ discrete values $(r_{\gamma})\in\Gamma^{\mathbb
{N}}$ independently drawn in $\Gamma$ with the probability density
$d\lambda(r)$. This sum hence resembles a Monte Carlo approximation
of the integral term, and we now show that our sums over populations
converge for almost all choices of $(r_{\gamma})\in\Gamma^{\mathbb
{N}}$ toward the integral, using an argument similar to the one we just
used to control $D_t(N)$. In detail, we show that $\mathcal{E}(\mathbb
{E}[\sup_{0\leq s\leq t}| E_s(N) |^2])$ converges toward $0$,
using the same method as that used for the convergence of the
mean-field term. Let us denote for the sake of compactness of notations
$F(s,r,r')$ the expectation $\mathbb{E}_Z[b(r,r',\bar
{X}_s^i(r),\break Z_{s-\tau(r,r')}(r')]$.

We have
\begin{eqnarray*}
&& \mathcal{E}\Bigl(\mathbb{E}\Bigl[\sup_{0\leq s\leq t}\bigl|
E_s(N) \bigr|^2\Bigr]\Bigr)
\\
&&\qquad \leq T \int
_0^t \mathcal{E} \Biggl(\mathbb{E} \Biggl[ \Biggl|
\int_{\Gamma} \frac{1}{P(N)} \sum_{\gamma=1}^{P(N)}
F(s,r_{\alpha
},r_{\gamma}) - \mathcal{E}_{r'} \bigl[F
\bigl(s,r_{\alpha},r'\bigr)\bigr] \Biggr| ^2\,ds
\Biggr] \Biggr).
\end{eqnarray*}
Similar to what was done for the term $D_t(N)$, since $\mathcal
{E}_{r'} [F(s,r_{\alpha},r')]$ is precisely the expectation of
$F(s,r_{\alpha},r_{\gamma})$ under the law of $r_{\gamma}$ over which
the sum is taken, developing the squared sum into a double sum over
populations (say, $\gamma$ and $\gamma'$), it is easy to show that,
because of the independence of the $r_{\gamma}$, that all terms that do
not correspond to $\gamma=\gamma'$, $\gamma=\alpha$ or $\gamma'=\alpha$
vanish, leaving less than $3P(N)$ possibly nonnull terms, and these
terms are uniformly bounded. Indeed, for $r_{\gamma }=r_{\gamma'}$ (the
case $r_{\gamma}=r_{\alpha}$ is treated in the same manner), we have
\begin{eqnarray*}
&& \mathcal{E}\bigl(\mathbb{E}\bigl[\bigl| F(s,r_{\alpha},r_{\gamma}) -
\mathcal {E}_{r'}\bigl[F\bigl(s,r_{\alpha},r'_{\gamma}
\bigr)\bigr] \bigr|^2 \bigr]\bigr)
\\
&&\qquad \leq 2 \mathcal {E}\bigl(\mathbb{E}
\bigl[\bigl| F(s,r_{\alpha},r_{\gamma}) \bigr|^2 +
\bigl|\mathcal{E}_{r'}\bigl[F\bigl(s,r_{\alpha},r'_{\gamma}
\bigr)\bigr] \bigr|^2 \bigr]\bigr)
\\
&&\qquad \leq2 \mathcal{E}\bigl(\mathbb{E}\bigl[\bigl| F(s,r_{\alpha},r_{\gamma})
\bigr|^2 + \mathcal{E}_{r'}\bigl[\bigl| F
\bigl(s,r_{\alpha},r'_{\gamma}\bigr)\bigr| ^2
\bigr] \bigr]\bigr)
\\
&&\qquad \leq 4 \tilde{K} \bigl(1+{C(s)}\bigr)
\end{eqnarray*}
implying eventually that
\begin{eqnarray*}
\mathcal{E}\Bigl(\mathbb{E}\Bigl[\sup_{0\leq s\leq t}\bigl|
E_s(N) \bigr|^2\Bigr]\Bigr) &\leq& \frac{4 T^2 \tilde{K}}{P(N)}
\bigl(1+{C(T)} \bigr).
\end{eqnarray*}
All together, we hence have
\begin{eqnarray*}
&& \mathbb{E}\Bigl[\sup_{0\leq s\leq t}\bigl| X^{i,N}_s-
\bar{X}^i_s(r_{\gamma }) \bigr|^2\Bigr]
\\
&&\qquad  \leq K'\int_0^{t} \max
_{j=1,\ldots, N} \mathbb{E}\Bigl[\sup_{-\tau\leq u \leq s}\bigl|
X^{j,N}_u-\bar{X}^j_u(r_{p(j)}) \bigr|^2\Bigr] \,ds
\\
&&\quad\qquad{}+C_1 \mathbbe (N) +\mathbb{E}\Bigl[\sup_{0\leq s\leq t}
\bigl| E_s(N) \bigr|^2\Bigr]
\end{eqnarray*}
valid for all $i\in\mathbb{N}$, and hence we have
\begin{eqnarray*}
&& \mathcal{E}\Bigl[\max_{i=1,\ldots,N}\mathbb{E}\Bigl[\sup _{0\leq
s\leq t}\bigl| X^{i,N}_s- \bar{X}^i_s(r_{\gamma}) \bigr|^2\Bigr] \Bigr]
\\
&&\qquad \leq K'\int_0^{t}\mathcal{E}
\Bigl[\max_{j=1,\ldots,N} \mathbb {E}\Bigl[\sup_{-\tau\leq u \leq s}
\bigl| X^{j,N}_u-\bar{X}^j_u(r_{p(j)})
\bigr|^2\Bigr] \Bigr] \,ds +C_1 \mathbbe (N) +
\frac{C_2}{P(N)},
\end{eqnarray*}
where $K'=4T (K_f^2+2L^2)$, $C_1=12 T^2 C$ and $C_2=16 T^2 \tilde
{K}(1+C(T))$ neither depend upon $N$ nor in the particular neuron
considered. By Gronwall's inequality, we hence obtain
\[
\mathcal{E} \Bigl[\max_{j=1,\ldots,N}\mathbb{E}\Bigl[\sup
_{-\tau\leq
s\leq t}\bigl| X^{j,N}_s-
\bar{X}^j_s(r_{\gamma}) \bigr|^2\Bigr]
\Bigr] \leq \biggl(C_1 \mathbbe (N) + \frac{C_2}{P(N)} \biggr)
\frac{e^{K'T}}{K'},
\]
which completes the proof.
\end{pf}

\begin{corollary}\label{corPropaChaosSpace}
Let $l\in\mathbb{N}^*$ and fix $l$ neurons $(i_1,\ldots,i_l) \in
\mathbb{N}^*$ with $i_{k}\neq i_{k'}$ for $k\neq k'$. Under the assumptions of
Theorem~\ref{thmPropagationChaosSpace}, the law of $(X^{i_1,N}_t,
\ldots, X^{i_l,N}_t, -\tau\leq t \leq T)$ converges toward
$m_t{(r_{p(i_1)})}\otimes\cdots\otimes m_t{(r_{p(i_l)})}$ for almost
all realization of the population locations $(r_{\alpha}, \alpha\in
\mathbb{N})$.
\end{corollary}

\begin{pf}
We have
\begin{eqnarray*}
&&\mathcal{E} \Bigl(\mathbb{E} \Bigl[ \sup_{-\tau\leq t \leq T} \bigl\llvert
\bigl(X^{i_1,N}_t, \ldots, X^{i_l,N}_t
\bigr) - \bigl(\bar{X}^{i_1}_t, \ldots, \bar{X}^{i_l}_t
\bigr)\bigr\rrvert ^2 \Bigr] \Bigr)
\\
&&\qquad \leq\sum_{k=1}^l\mathcal{E}
\Bigl(\mathbb{E} \Bigl[ \sup_{-\tau\leq t \leq T} \bigl\llvert
X^{i_k,N}_t-\bar{X}^{i_k}_t\bigr\rrvert
^2 \Bigr] \Bigr)
\leq l \biggl(C_1 \mathbbe (N) + \frac{C_2}{P(N)}
\biggr) \frac{e^{K'T}}{K'},
\end{eqnarray*}
which tends to zero as $N$ goes to infinity; hence the law of
$(X^{i_1,N}_t, \ldots, X^{i_l,N}_t,\break  -\tau\leq t \leq T)$ converges
toward that of $(\bar{X}^{i_1}_t, \ldots, \bar{X}^{i_l}_t, -\tau\leq t
\leq T)$ whose law is equal by definition to
$m(t,r_{p(i_1)})\otimes\cdots\otimes m(t,r_{p(i_l)})$. Since the
expectation $\mathcal{E}$ of the distance between the processes
considered tend to zero, these processes\vadjust{\eject} converge $\mathbb{P}'$-almost
surely, that is, for almost all realizations of the space locations,
which ends the proof.
\end{pf}

\begin{Imremark*}
The speed of convergence toward the mean-field equation is hence
governed by $\mathbbe (N)$ and $1/P(N)$. In the case of a finite
number\vadjust{\goodbreak} of populations, the speed of convergence is hence driven by the
size of the smallest population. In the infinite population case, the
speed of convergence toward the mean-field limit is a balance between
the averaged number of neurons in each population through the term
$\mathbbe (N)$, and the total number of populations through the term
$1/P(N)$. The first term quantifies the speed at which averaging
effects occur in the network and is related to the averaged inverse
number of neurons in each population. The other term controls the
convergence of the interaction related to all populations toward an
effective interaction term given by an integral over $\Gamma$ of
mean-field interactions, that is, convergence of finite-populations
networks toward their continuous limit. For networks with homogeneous
population sizes, $\mathbbe (N)$ will be approximately equal to
$P(N)/N$. The optimal network size ensuring the fastest convergence in
that case hence corresponds to $P(N)\sim\sqrt{N}$ [minimizing the
functional $x\mapsto P(x)/x + 1/P(x)$], and in that case the
convergence will be in $1/\sqrt{N}$, and we conjecture that this speed
of convergence is optimal (though we did not achieve to prove it). This
convergence speed is hence very slow compared to finite-size networks
and usual mean-field limits in which the speed of convergence is of
order $1/N$.
\end{Imremark*}

\section{Neural fields equations in action}\label{secFiringRates}
It is folklore that McKean--Vlasov limits have dynamics that are
complex to analyze. Very refined methods are generally set up to
analyze the behavior of the system in the mean-field limit, such as
entropy methods or spectral methods; see, for example,
\cite{villani02}. This statement could be even more true in our
spatialized context, and the present, general approach might appear to
be bounded to remain formal.

Fortunately, for relevant neuroscience applications, it happens that
solutions to these equations are not out of reach. This is the topic of
a companion article~\cite{touboulNeuralFieldsDynamics11} where networks
of firing-rate neurons (see Appendix~\ref{appendFiringRates}), are
considered, the neuronal model usually considered for neural fields
analysis. Let us briefly review here the main results of that article
and concretely use the proposed approach to analyze the dynamics of a
simple network.

Considering firing-rate neurons, we show
in~\cite{touboulNeuralFieldsDynamics11} that the solutions of the
mean-field equations are Gaussian processes when the initial condition
is as well (and equilibria are Gaussian) and that their mean $M(r,t)$
and standard deviation $v(r,t)$ (fully describing the process since the
covariance is a simple function of these two quantities in that case)
reduces to the set of deterministic delayed integro-differential
equations
\begin{equation}
\label{eqMomentsGaussian}
\cases{ \displaystyle\partial_t M(r,t) =-\frac{1}{\theta(r)}
M(r,t)+I(r,t)
\vspace*{5pt}\cr
\displaystyle\qquad\hspace*{30pt}{} +\int_{\Gamma}J\bigl(r,r'
\bigr)
\vspace*{5pt}\cr
\hspace*{75pt}{}\times  F\bigl(r',M\bigl(r',t-\tau\bigl(r,r'
\bigr)\bigr),v\bigl(r',t-\tau \bigl(r,r'\bigr)\bigr)
\bigr)\lambda\bigl(r'\bigr)\,dr',
\vspace*{5pt}\cr
\displaystyle\partial_t v(r,t) =-\frac{2}{\theta(r)} v(r,t)+\sigma(r)^2,}\hspace*{-30pt}
\end{equation}
where $F(r,x,y)$ denote the expectation of $S(r,U)$ for $U$ a Gaussian
random variable of mean $x$ and variance $y$, and can be made explicit
for particular choices of sigmoids $S$. These equations are consistent
with the heuristically derived extremely widely used Wilson--Cowan
models for finite-populations neural
assemblies \mbox{\cite{wilson-cowan72,wilson-cowan73}} in the limit where
noise levels vanish. These equations are shown to be well-posed, and
grant access to the dynamics of the network.
In~\cite{touboulNeuralFieldsDynamics11}, the choice of the parameters,
driven by biological constraints, did not reveal any qualitative effect
of the delays on the solutions except during transient phases.

In order to illustrate how the use of the present approach can be used
to uncover the dynamics of the neural field, we proceed to the analysis
of a single population network with inhibitory interactions (i.e.,
negative interactions), a case that was not treated
in~\cite{touboulNeuralFieldsDynamics11} and which will turn out show a
particularly rich variety of behaviors as a function of delays.

To this end, let us fix the parameters of the system. We consider
$\Gamma=\mathbb{S}^1$ the 1-dimensional torus, and $\lambda$ the
uniform distribution on it. We consider that $S(r,x)=\int_0^{gx}
e^{-x^2/2}/\sqrt{2\pi}=:\operatorname{erf}(gx)$, $\theta(r)=1$ and
$\sigma$ independent of $r$, and one can easily show by changing
variables that $F(x,y)=\operatorname{erf}(gx/\sqrt{1+g^2y})$. We
further fix $J(r,r')=\bar{J} e^{-| r-r'|/\delta}$ ($\delta$
represents the typical connectivity length in the neural field) and
$\tau(r,r')=| r-r' |/ c +\tau_s$ ($c$ represents the speed of
transmission in the neural field and $\tau_s$ the typical transmission
time of the synapse).

Since $F(0,y)=0$ for any $y\in\mathbb{R}$, the Gaussian solutions with
zero mean and standard deviation $\sigma^2/2$ are stationary solutions
of the system that are spatially homogeneous in law (i.e., their law
does not depend on the space variable). Characterizing the stability of
this solution consists of analyzing the characteristic roots equation
of the linearized system around the spatially homogeneous stationary
solution. Computing the eigenvalues of the integral convolution
operator similarly to~\cite{touboulNeuralFieldsDynamics11},
Section~3.1, we obtain the \emph{dispersion relationship}
\[
\xi+1= F_0' \frac{e^{-\xi\tau_s}(1-e^{-(1/\delta+ \xi/
c)})}{1/\delta+ \nu/ c+\mathbf{i}2\pi k }
\]
for $k\in\mathbb{Z}$ and $F'_0=\frac{g}{\sqrt{1+g^2v_0}}
\frac{1}{\sqrt{2\pi}}$. The spatially homogeneous equilibrium is stable
if and only if all solutions $\xi$ to the dispersion relationship
(characteristic roots) have negative real parts. A Turing bifurcation
point is defined by the fact that there exists an integer $k$ such that
$\Re(\xi)=0$. It is said to be \emph{static} if at this point
$\Im(\xi)=0$, and \emph{dynamic} if $\Im(\xi)=\omega_k\neq0$. In that
latter case, the instability is called a Turing--Hopf bifurcation, and
generates a global pattern with wavenumber $k$ moving coherently at
speed $\omega_k/k$ as a periodic wavetrain.

\begin{figure}

\includegraphics{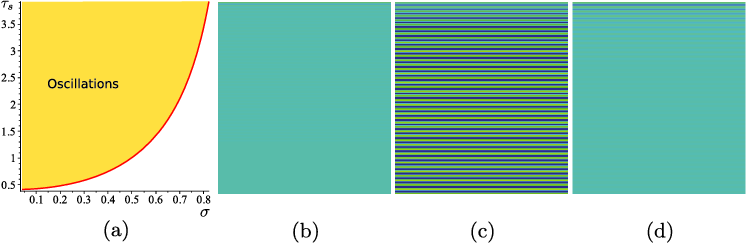}

\caption{Turing--Hopf bifurcations and delay-induced synchronization,
$J=-3$, $g=3$, $\delta=c=1$. \textup{(a)}:~bifurcation diagram, shows a
transition from stationary to periodic activity as delays are increased
\textup{(b)${}\to{}$(c)}.
\textup{(b)}~$\sigma=0.1$, $\tau_s=0.4$,
\textup{(c)}~$\sigma=0.1$, $\tau_s=0.5$,
\textup{(d)}~$\sigma=0.3$, $\tau_s=0.5$.
When noise is increased, synchronization is lost \textup{(c)${}\to
{}$(d)}. \textup{(b)--(d)}: spatio-temporal dynamics as a function of space
(abscissa) and time (ordinate).}
\label{figDelays}
\end{figure}

\begin{figure}

\includegraphics{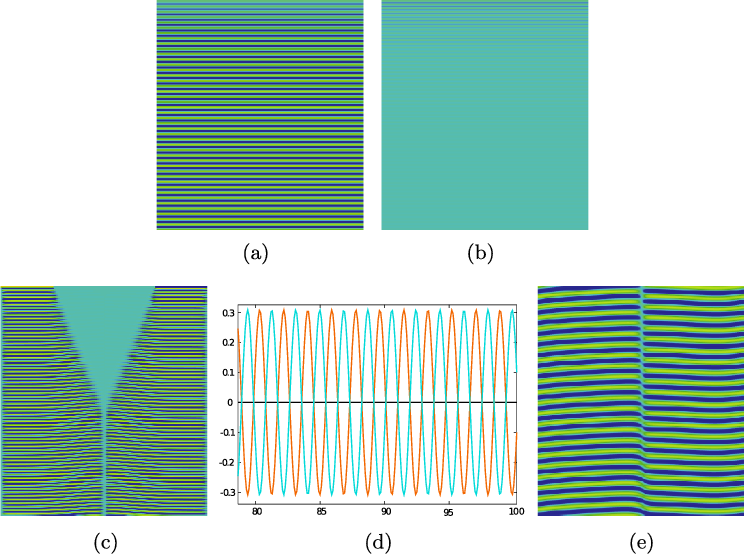}

\caption{\textup{(a)}~$\delta=1$,
\textup{(b)}~$\delta=5$,
\textup{(c)}~and \textup{(d)}~antiphase,
\textup{(e)}~synchronization.
Spatial effects: increasing $\delta$ destroys the
synchronization. For $\delta=1$
[case~\textup{(c)} of Figure~\protect\ref{figDelays}], choosing nonspatially
homogeneous initial
conditions yields to complex situations, for instance, an antiphase
synchronization during
long transients (bottom row).
\textup{(c)}: $t\in[0,200]$,
\textup{(d)}: orange: $M(t,0.1)$, blue: $M(t,0.9)$, black: $M(t,0.5)$,
\textup{(e)}: $t\in[600, 650]$. The synchronization becomes visually perfect
for times above 1500.}\label{figDelaysSpace}
\end{figure}

Possible Turing--Hopf bifurcations hence arise when there exists
$\omega_k>0$ such that
\[
\mathbf{i}\omega_k+1= F_0'
e^{-\mathbf{i}\omega_k\tau_s}Z_k(\omega)
\]
with $Z_k(\omega)=\frac{(1-e^{-(1/\delta+ \xi/
c)})}{1/\delta+ \nu/ c+\mathbf{i}2\pi k }$, which yields
bifurcation curves (parametrized by $\omega$) in the parameter space
\[
\cases{ \sigma^2 =\displaystyle{\frac{2}{g^2} \biggl(-1+
\frac{\bar
{J}^2g^2| Z_k(\omega)|^2}{2\pi(1+\omega^2)} \biggr)},
\vspace*{5pt}\cr
\tau_s = \displaystyle{
\frac{-\arctan(\omega) + \operatorname{Arg}(F_0' Z_k(\omega)) + 2m\pi}{\omega}}.}
\]
This provides a curve of Turing--Hopf bifurcations corresponding to
transitions from stationary independent solutions to perfectly
synchronized independent solutions, as displayed in
Figures~\ref{figDelays} and~\ref{figDelaysSpace}.
%
In Figure~\ref{figDelays}, we display the bifurcation curve in the
parameter space $(\sigma,\tau_s)$ for a specific set of parameters.
This curve has a convex shape. Small enough delays hence correspond to
stationary solutions. Increasing delays yields periodic activity, which
disappears as noise is increased. This example shows the importance of
delays in the qualitative dynamics of the neural field. The typical
connectivity length also shapes the qualitative dynamics of the neural
field, as shown in Figure~\ref{figDelaysSpace}. This variety of
behaviors correspond to bifurcations corresponding to a wavenumber
$k=0$, and correspond to spatially homogeneous solutions. Nontrivial
spatial structures can be searched for considering nonspatially
homogeneous initial conditions. In this case, a number of
\mbox{complex}
spatio-temporal behaviors can appear, such as the metastable
polychronization shown in Figure~\ref{figDelaysSpace}, where the neural
field splits into two clusters oscillating in antiphase during very
long transient periods before a sudden synchronization of the whole
neural field.

%

\section{Discussion}\label{secdiscussion}
In this paper, we addressed the problem of the asymptotic behavior of
networks composed of a large number of neuronal assemblies in a
particular asymptotic regime, the \emph{neural-field limit}. We took
into account\vadjust{\goodbreak} a number of specificities relevant to neuronal dynamics:
intrinsic noise at the level of each neuron, the spatial structure and
propagation delays. We demonstrated that for \mbox{a~relatively} general class
of models, that includes the most prominent models in neuroscience
(reviewed in Appendix~\ref{secMotiv}), the propagation of chaos
property took place and showed convergence of the mean-field equations
toward mean-field equations of a new type, analogous to the classical
McKean--Vlasov equations, but including delayed interactions, a spatial
integration term and a singular spatio-temporal stochastic process, the
spatially chaotic Brownian motion.

The question of the scale at which relevant phenomena occur is
essential to the modeler. Descriptions coarser than our neural field
limit, for instance, those involving finite numbers of populations,
correspond to cases where our measure $\lambda$ is a sum of Dirac
masses. This case can be seen as a particular case of the present
analysis, and hence the propagation of chaos occurs and network
equations converge toward mean-field equations that correspond to a
finite system of delayed McKean--Vlasov equations. In contrast, scales
finer than the neural field limit (taking, e.g., into account possible
individual heterogeneities between neurons) are not covered by the
analysis and seem relatively hard to understand. It is likely that the
dynamics of such networks will be considerably distinct from that of
networks in the neural field regime. The neural-field regime seems
particularly well suited to describe the activity of large neuronal
assemblies, since it was observed that population sizes are orders of
magnitude larger than the total number of
populations~\cite{changeux06}. Moreover, it seems to be at the scale of
biological recordings and phenomena such as the emergence of patterns
of activity in the cortex. We illustrated how such an analysis could be
rigorously developed with a simple example in
Section~\ref{secFiringRates}. More relevant states may be analyzed with
this model, since the usual heuristic equations that were successfully
used in a number of situations~\cite{ermentrout-cowan79} are
compatible, in the zero noise limit, with our equations, and the
rigorously derived model will shed new light on the role of noise in
such neuronal systems, but also on the individual behaviors of neurons.
For instance, the propagation of chaos property ensures that finite
sets of neurons are independent in the neural field limit. This result
contradicts the classical view considering that since neurons of the
same population are highly connected and receive similar input, their
activity shall be correlated. However, with recent experimental
findings using high-quality
recordings~\cite{ecker-berens-etal10,renart-de-la-rocha-etal10} showed
that levels of correlations between two neurons (of the same population
or not) were extremely small, way below what was usually considered.
The propagation of chaos hence offers a universal explanation to this
phenomenon.

A number of open questions remain widely open in the theoretical
understanding of the behavior of neural fields and large-scale neural
networks. For instance, a particularly\vadjust{\goodbreak} interesting phenomenon is the
plasticity of neuronal connections. Considered constant and homogeneous
in the present manuscript, it happens that the synaptic coefficients
describing pairwise interactions between neurons evolve, very slowly,
as a function of the network activity and in particular as a function
of the correlations between the activity of pairs of neurons. This kind
of phenomena was never considered in the mathematical literature, and
seems relatively rich. In particular, this mechanism can break the
propagation of chaos property and yield weakly correlated states. This
is a problem we are currently investigating.

\begin{appendix}
\section{Neuron models}\label{secMotiv}
For the sake of completeness we quickly review in this
appendix different classical neuron models motivating the present study. This
appendix takes a mathematical viewpoint, is obviously very selective
and lacunar. The interested reader will find more details in classical
neuroscience textbooks, for
example,~\cite{kandel-schwartz-etal00,ermentrout-terman10}. Basically,
neurons are electrically excitable cells whose activity, measured
through the voltage of the cell (difference of electrical potential
between the intracellular and extracellular domains), is governed by
ionic transfers through specific proteins (ion voltage-gated ion
channels) located on the cellular membrane. We present here detailed
neuron models (Appendix~\ref{appendcondBased}) that approximate the
biophysics of ion channels, and firing-rate models
(Appendix~\ref{appendFiringRates}) that reproduce qualitatively the
dynamics of the firing rate of neurons and that are used in the
application Section~\ref{secFiringRates}.

\subsection{Hodgkin--Huxley and Fitzhugh--Nagumo models}\label{appendcondBased}
Probably the most biologically relevant, versatile
and precise neuron model is the Hodgkin--Huxley (HH)
model~\cite{hodgkin-huxley52}. This model describes the membrane
potential $v$ of a~neurons as a function of the dynamics of several
ionic currents that enter or exit the cells through voltage-gated
channels. The mathematical description we choose here involves Langevin
approximation of the random proportion open of ion channel; see, for
example,~\cite{goldwyn-shea-brown11} and references therein. The
proportion of open channels satisfies in that model a stochastic
differential equation,
\[
dx_t=\bigl(A_x(v) (1-x)-B_x(v) x\bigr) \,dt
+ \sqrt{A_x(v) (1-x)+B_x(v) x}\chi (x)
\,dW^x_t,
\]
where $W^x_t$ are independent standard Brownian motions, $A_x(v)$ and
$B_x(v)$ are smooth bounded functions accounting, respectively, for the
opening and closing probability intensity of a given channel and
$\chi(x)$ is a function vanishing outside $[0,1]$ to ensure that the
variables $x$ remain in $[0,1]$ (since these variables describe
proportions). Generally, three ionic currents (and channels) are
considered: potassium ($m$), calcium activation ($n$) and inactivation
($h$) and Ohmic leak current, $I_L$ (carried by Cl$^-$ ions).
Considering that the neuron receives an external current composed of a
deterministic part $I(t)$ and a white noise with standard deviation
$\sigma_{\mathrm{ext}}$, the voltage is governed by the equation
\begin{equation}
\label{eqHodgkinHuxley}
\cases{
\displaystyle C \,d{v}_t = \bigl(I(t) - \bar{g}_K
n^4 (v-E_K) - \bar{g}_{Na} m^3
h (v-E_{Na}) - \bar{g}_L (v-E_L) \bigr) \,dt
\vspace*{3pt}\cr
\hspace*{36pt}{} +\sigma_{\mathrm{ext}} \,dW_t,
\vspace*{3pt}\cr
\displaystyle d{x}_t =
\bigl(A_x(v) (1-x) -B_x(v) h \bigr) \,dt +
\sigma_x(v,n) \,dW^x_t x\in\{n, m, h\}.}\hspace*{-30pt}
\end{equation}
This model satisfies assumptions~\hyperref[AssumpLocLipschSpace]{(H1)}
and~\hyperref[AssumpLinearGrowth]{(H4)} used in the general theory,
since though polynomial nonlinearities arise in the dynamics, the
boundedness of the variables $(n,m,h)$ ensure Lipschitz continuity and
linear growth. Assumption~\hyperref[AssumpSpaceContinuity]{(H5)} is not
satisfied since the noise depends on the state of the neuron. This
refinement does not make the proofs substantially more intricate as
discussed in Appendix~\ref{appendStoSynapses}.

The HH model is often too complex for practical purposes, and several
reductions were proposed. A particularly interesting one is the
Fitzhugh--Nagumo (FN) bidimensional model \cite{fitzhugh69} capturing
from the biological viewpoint the most prominent behaviors of the
Hodgkin--Huxley model. From the mathematical viewpoint, it is important
to specify this model since that model does not satisfy
assumptions~\hyperref[AssumpLocLipschSpace]{(H1)}
and~\hyperref[AssumpLinearGrowth]{(H4)}, and motivates the additional
mathematical developments of Appendix~\ref{appendStoSynapses}. This
model describes the evolution of the membrane potential variable $v$
and a slower recovery variable $w$, through the equations
\begin{equation}
\label{eqFNStoch} \cases{ dv_t=\bigl(P(v_t) -
w_t+I\bigr) \,dt + \sigma_v\,dW^v_t,
\vspace*{3pt}\cr
dw_t=a (b v_t-w_t) \,dt +
\sigma_w \,dW^w_t,}
\end{equation}
where $P(v)=v(1-v)(v-a)$, generally chosen $f(v)=v-v^3$.

The state of the neuron $X$ in our abstract
model~(\ref{eqNetworkSpace}) in the HH model is given by $(v,n,m,h)$
and for the FN model by $(v,w)$, and their intrinsic dynamics is
enclosed in the functions $f$ and $g$.

The communication between neurons is maintained by two possible types
of synapses: electrical or chemical. Electrical synapses, in charge of
rapid and stereotype signal transmission, operate through direct
contact of the intracellular domain of the two communicating cells
through specialized protein structures called gap-junctions. The ions
passively flow from one neuron to the other: the interaction is not
delayed, and the current produced by neuron $j$ on neuron $i$ is equal
to $J_{ij}(v^j_t-v^i_t)$ where $J_{ij}$ is called the synaptic
conductance [this defines our interaction function $b$ in the abstract
model~(\ref{eqNetworkSpace})]. When including the dependence on $v^i_t$
in the drift function, the interaction function $\sum_j J_{ij}v^j_t$
clearly satisfies assumptions~\hyperref[AssumpLocLipschbSpace]{(H2)}
and~\hyperref[AssumpbBoundSpace]{(H3)},
and~\hyperref[AssumpSpaceContinuity]{(H5)} as soon as the dependence of
$J_{ij}$ with respect to space is sufficiently regular. The chemical
synapse is the most common type of interconnection. When a spike is
fired from a pre-synaptic neuron $j$, it is transported through the
axons to the synaptic button where it is transmitted to neuron $i$
through a complex process of release of neurotransmitter (from $j$)
binding to specific receptors on neuron $i$. The transmission takes a
time $\tau_{ij}$ in the order of a~few milliseconds. Similar to HH
ion channels dynamics, the proportion of open neurotransmitter channels
$y^i$ has the dynamics (see~\cite{destexhe-mainen-etal94b}),
\[
dy^j_t= \bigl(A S\bigl(v^j\bigr)
\bigl(1-y^j(t)\bigr)-D y^j(t) \bigr) \,dt + \sigma
_Y\bigl(v^j,y^j\bigr) \,dW^{j,y}_t
\]
with $S$ is a smooth sigmoidal function. In our abstract model, the
variable $y^i$ is added to the state $X^i$ of neuron $i$, and the
functions $f$ and $g$ take into account that dynamics. The synaptic
current induced at time $t$ on neuron $i$ by the arrival of a~spike
from neuron $j$ (fired at time $t-\tau_{ij}$) is equal to $J_{ij}
y^j(t-\tau_{ij}) (v^i(t)-v_{\mathrm{rev}})$ governing our interaction
function $b$ clearly satisfying
assumptions~\hyperref[AssumpLocLipschbSpace]{(H2)}
and~\hyperref[AssumpbBoundSpace]{(H3)},
and~\hyperref[AssumpSpaceContinuity]{(H5)} as soon as the dependence of
$J_{ij}$ with respect to space is sufficiently regular.

The synaptic efficacies $J_{ij}$ of electrical or chemical synapses are
given by the connectivity of the cells. Such functions are generally
considered continuous functions $J(r_i,r_j)$ depending on the
population of $i$ and $j$.

Putting all these elements together and assuming that all the
parameters of the equations only depend on the neural populations of
the cells involved, we can write the equation of a network of FN
neurons with chemical synapses, external and synaptic noise,
\begin{equation}
\label{eqFNNetwork}
\cases{ \displaystyle dv^i_t = \Biggl(P
\bigl(v^i_t\bigr)+I^{i}(t)+ \sum
_{j=1, j\neq i}^N \bigl(J_{ij}
y^{j}(t-\tau_{ij}) \bigl(v^i_{t}-v_{\mathrm{rev}}
\bigr) \bigr) \Biggr)\,dt,
\vspace*{3pt}\cr
\displaystyle d{w}^i_t = a_{\alpha}
\bigl(b_{\alpha} v^i_t-w^i_t
\bigr)\,dt,
\vspace*{3pt}\cr
\displaystyle dy^i_t = \bigl(A_{\alpha}
S_{\alpha}\bigl(v^i_t\bigr) \bigl(1-y^i_t
\bigr)-D_{\alpha} y^i_t\bigr) \,dt+
\sigma_Y(v,y) \,dB^{i, Y}_t.}
\end{equation}
A similar (but more complex) expression is obtained for the HH model
using equations~(\ref{eqHodgkinHuxley}) and with distributed delays.

\subsection{Stochastic firing-rates models}\label{appendFiringRates}
A phenomenological neuron model consists of considering that neurons
interact through their mean firing-rate. The firing-rate model
considers that the membrane potential has a linear dynamics, and its
mean-firing rate is a smooth sigmoidal transform of the membrane
potential $S(r_{\alpha},\cdot)$ depending on the neural population
$\alpha$. In other words, an incoming firing rate provokes postsynaptic
potentials that linearly sum. The neurons receive additional inputs
that are the sum of a deterministic current $I(r_{\alpha},t)$ and noise
$\sigma(r_{\alpha})\,d{W}^i_t$. The network equations hence read
\begin{eqnarray*}
dV^{i}(t) &=& \Biggl( -\frac{1} {\theta(r_{\alpha})} V^{i}(t) +
I(r_{\alpha},t)
\\
&&\hspace*{6pt}{} + \sum_{\gamma=1}^{P}
J_{\alpha\gamma} \frac
{1}{N_{\gamma}} \sum_{j, p(j)=\gamma} S
\bigl(r_{\gamma},V^j(t-\tau _{\alpha\gamma})\bigr) \Biggr) \,dt
  + \sigma(r_{\alpha})\,dW^i_t.
\end{eqnarray*}

It is easy to check that
assumptions~\hyperref[AssumpLocLipschSpace]{(H1)}--\hyperref[AssumpSpaceContinuity]{(H5)}
are satisfied for the firing-rate model.

\section{Generalized models}\label{appendStoSynapses}
In the main section we choose to concentrate on the cornerstone
mathematical problems arising in the modeling of neural fields, and
choose to deal with relatively general models, yet simplified. Indeed,
as discussed in Appendix~\ref{secMotiv}, we see that two technicalities
were not taken into account in our general analysis. These were (i)
nonglobally Lipschitz drift that do not satisfy the linear growth
condition (for Fitzhugh--Nagumo models) and (ii) state-dependent
diffusion coefficients. 

{\bf Update} While the latter point can be easily taken into account with
usual methods, the question of nonglobally Lipschitz continuous drifts and 
diffusions functions are relatively delicate and non-classical. As noted 
in~\cite{BFT}, the result presented to be potentially true in the original version 
would be incorrect under the too weak assumptions initially made, and counter-examples 
exist~\cite{scheutzow}. In particular, the localization and truncation method suggested
is not correct. 

We add that the only model for neuroscience applications would motivate to extend the 
results to locally-Lipschitz drifts is the Fitzhugh-Nagumo network equation~(\ref{eqFNNetwork}).
The reader interested in extending the results to the FitzHugh-Nagumo model is referred to~\cite{BFT}
where a probabilistic method is proposed in order to show existence and uniqueness of solutions 
on a finite time interval in a non spatially extended setting. The methods could possibly extend to 
the present setting. 
An alternative approach based on functional analysis methods~\cite{quininao} developed also in
a non-spatialized setting may also generalize to the present case to the Fitzhugh-Nagumo setting, 
and could allow proving global in time results of existence and uniqueness of solutions, existence 
of stationary solutions and their stability in the weak coupling regime. 

\end{appendix}

\section*{Acknowledgments} The author warmly acknowledges Cl\'ement
Mouhot for very interesting discussions on the links between this
research and kinetic theory. A~special acknowledgment goes to Philippe
Robert for great discussions and advice on the content and writing of
the paper. He also acknowledges Gerard Ben Arous and Marc Yor for
interesting discussions and suggestions.




\printaddresses

\end{document}